\documentclass[review,letterpaper]{elsarticle}

\usepackage[english]{babel}
\usepackage[utf8]{inputenc}

\usepackage{geometry}
\usepackage{multicol,float}
\usepackage{subcaption}
\usepackage{color}
\usepackage{xcolor,colortbl}
\usepackage{multirow}
\usepackage{array}
\usepackage{amsmath}

\usepackage{amsfonts}
\usepackage{amssymb}
\usepackage{graphicx}
\usepackage[ruled, vlined, linesnumbered]{algorithm2e}
\usepackage{mathtools}
\usepackage{soul}
\usepackage{cancel}
\usepackage[normalem]{ulem}
\usepackage{booktabs}
\usepackage{multirow}
\usepackage{url}
\usepackage{tabularx}
\usepackage{hyperref}

\colorlet{shadecolor}{gray!20}

\title{Optimizing Accessibility and Maximum Travel Times in Transit Networks: The Hub Covering Location Problem}
\author[uca]{Carmen-Ana Domínguez-Bravo}
\author[uca]{Elena Fernández}
\author[unab,catlec]{Armin Lüer-Villagra}
\address[uca]{Statistics and Operations Research Department, Universidad de Cádiz, Puerto Real, Spain}
\address[unab]{Department of Engineering Sciences, Universidad Andres Bello, Talcahuano, Chile}
\address[catlec]{Center for Advanced Transportation, Logistics, and Economic Competitiveness (CATLEC), Chile}
\date{}

\begin{document}

\begin{abstract}
This work introduces an optimization model to support decision-making in the design of integrated transit systems, with the aim of improving urban mobility.
We assume that daily user trips consist of two types of segments: short segments where users apply their own transport between their homes or workplaces and the public transit access points, and a longer segment on the transportation network.\\
Users' accessibility is ensured by imposing an upper limit on the travel times of the short segments. We assume that only users within a given coverage radius will employ the system. On the other hand, system effectiveness is ensured by limiting travel time on the transportation network.\\
Two mixed-integer linear programming formulations and a matheuristic procedure are proposed, and their efficiency is evaluated through computational experiments. The results allow us to study the structure of the solutions obtained and analyze the effects of the parameters considered.
\end{abstract}


\maketitle

\section{Introduction}\label{sec:intro}
Since the last decades of the last century, the size of large cities has expanded seemingly without limit. Daily commuting between residential areas and workplaces has created a vicious cycle: as traffic congestion increases, travel times become longer and less predictable. In response, more people turn to private vehicles in an attempt to regain control over their journeys--further aggravating congestion.

In 2019, people in employment in the EU had an average commuting time of 25 minutes\footnote{\url{http://ec.europa.eu}}. 
According to the US Census Bureau\footnote{\url{https://www.census.gov/}}, the average one-way travel time in 2024 was 27.2 minutes, higher than the 26.8 minutes in 2023. The same source indicates that the percentage of workers with a one-way travel time of 60 minutes or more in 2024 was 9.3 percent, higher than the 8.9 percent in 2023. Global estimates by the Worldmapper Archive\footnote{\url{http://archive.worldmapper.org/posters/worldmapper_map141_ver5.pdf}} indicate that, in 2005, the world average commuting time was 40 minutes, one-way.

It is therefore unsurprising that urban mobility has become a major concern for city planners. Key considerations include improving the integration of public transportation systems by enhancing accessibility and reliability, providing greater control over commute times, and ensuring adequate physical spacing and comfort conditions to promote user well-being \cite{30-mins}. An overview of public transit can be found in \cite{DesaulniersHickman}. Among the relevant models proposed to enhance the accessibility of transit systems, we can mention location/set covering models (see, e.g., \cite{Murray2003}).

In this paper, we introduce a new optimization model to support decision-making in the design of integrated transit systems that improve accessibility while limiting user travel times. For the users of the proposed systems, we assume that daily trips consist of short segments made by their own means of transport (e.g., walking, cycling, or similar) between their homes or workplaces and the public transit system, combined with a longer segment within the public transportation network that connects the network entry and exit points. 
Accessibility is ensured by imposing an upper limit on the travel distance between the endpoints of the short segments. That is, we assume that potential users located beyond this threshold, i.e., radius, will not employ the public transportation system. Conversely, we also assume that users within the access coverage radius will benefit from the system, whose effectiveness is guaranteed because the travel times of the longer segments will remain below a prescribed maximum. Decision-makers must therefore strike an appropriate balance between the installation costs of such a system and the number of users it can effectively serve. This model will be referred to as the Hub Covering Location Problem (HCLP). 

The HCLP addresses strategic decisions related to the design of a transit network, integrating network design and covering decisions. Indeed, networks associated with low activation costs are unlikely to provide adequate accessibility for many users, whereas increasing accessibility necessarily entails higher activation costs. In our approach we penalize non-acessible areas, so the HCLP aims at finding a trade-off between accessibility and installation costs by minimizing the sum of the total activation costs plus the total penalties. Note that minimizing the penalties for uncovered areas is equivalent to maximizing the captured demand.

Even though the primary applications of the HCLP lie in ground transportation systems, it also applies to air transportation, particularly in the design of air transport networks where user itineraries combine short domestic origin–destination legs with longer flights (or combinations of flights) connecting them.

As we will see, the HCLP has some similarities with classical Hub Location Problems (HLPs) and, in particular with hub covering problems~\citep{Campbell1994}. However, a key difference, which clearly differentiates it from this class of problems, is the modeling of demand. In HLPs, demand is associated with commodities, i.e., with OD pairs, whereas in the HCLP, demand is associated with individual nodes. That is, in the HCLP, a node is served, i.e. captured by the network, if it is within the coverage radius $\delta$ of some activated facility and potential destinations of served nodes are modeled only implicitly. In particular, in the HCLP, any pair of served nodes determines a served commodity. This way of handling demand is also reflected in the implicit modeling of routing costs. Due to economies of scale, in hub-and-spoke networks, the routing costs of flows through interhub arcs are typically discounted, leading to lower (average) interhub routing costs. On the contrary, the objective function of the HCLP  does not explicitly consider the routing costs through interhub arcs, which are only implicitly taken into account by limiting the travel times within any activated facilities to a maximum value, $T_{max}$.

The main contributions of this paper are the following:
\begin{itemize}
\item We introduce the HCLP, a novel problem involving two simultaneous covering considerations: network access distance and maximum network transit time.

\item We study optimality conditions of the  HCLP and present two alternative mixed-integer linear programming (MILP) formulations. The first one, $F_1$, uses classical four-index path variables for the arcs in the path connecting any given pair of hubs.  The second one, $F_2$, uses three-index variables to model a spanning tree rooted at each of the hubs of the transportation system.

\item We moreover propose a two-phase solution algorithm. It is a matheuristic, which in each iteration solves a simplified HCLP consisting of two independent subproblems: a \textit{Location and Covering} (LC) subproblem and a  \textit{Network Design} (ND) subproblem. The LC phase determines a set of hubs to activate together with a set of users to be captured by the transportation system, and can be seen as a variable neighborhood search. For this, the LC formulation includes no-good cuts that prevent the same set of hubs from being activated across different iterations, plus additional constraints that restrict the search to hubs located within a certain radius of those already activated. The ND phase builds a transit network using the hubs selected in the first phase. Both phases are optimally solved through \textit{ad hoc} MILP formulations derived from $F_2$ , which, as we will see, outperforms $F_1$. 

\item Extensive computational experiments have been run on data sets adapted from the literature  to compare the two formulations and to analyze the performance of the solution algorithm. Furthermore, we develop sensitivity analyses to examine the effect of the different parameters on the solutions structure and to gain insight on the proposal.
\end{itemize}

The remainder of this paper is structured as follows. Section~\ref{sec:literature} presents a review of the related literature. Section \ref{sec:prob-def} provides the problem definition. Then, Sections \ref{sec:4i} and \ref{sec:3i} provide four- and three-index formulations, respectively. After that, Section \ref{sec:MathH} describes the two-phase solution algorithm for the problem. Section \ref{sec:compu} provides extensive computational experience. Finally, Section \ref{sec:conclu} presents the conclusions and future research directions.

\section{Literature review}\label{sec:literature}

The HCLP comprises simultaneous locational and network design decisions. Thus, it can be framed as a general network design problem \citep{Contreras2012}. Table~\ref{tab:literature} shows the classes of problems closest to the HCLP, which are addressed in the remainder of this section.

\begin{table}[H]
\caption{Comparison of our problem (HCLP) and related problem classes and subclasses.}
\label{tab:literature}
\resizebox{\textwidth}{!}{%
\begin{tabular}{llccccccc}
\textbf{Class}                     & \textbf{Subclass} & \multicolumn{1}{l}{\textbf{Location}} & \multicolumn{1}{l}{\textbf{Net. Design}} & \multicolumn{1}{l}{\textbf{Fixed costs}} & \multicolumn{1}{l}{\textbf{Variable costs}} & \multicolumn{1}{l}{\textbf{Spokes}} & \multicolumn{1}{l}{\textbf{Interfacility}} & \multicolumn{1}{l}{\textbf{Maximum time}} \\
\toprule
\multirow{1}{*}{Facility location} & Covering        & \checkmark &             & \checkmark  &            &            &            & \checkmark          \\
\midrule
\multirow{2}{*}{Hub location}      & Incomplete      & \checkmark & \checkmark  &             &  \checkmark          & \checkmark & \checkmark &                     \\
                                   & Hub covering    & \checkmark & \checkmark  & \checkmark  &            & \checkmark & \checkmark & \checkmark            \\
\midrule
\multirow{1}{*}{Network design}    & Conn. Fac. Loc. & \checkmark & \checkmark  & \checkmark  &             &            & \checkmark   &           \\                  
\midrule
This work                          & HCLP           & \checkmark  & \checkmark   & \checkmark &            &            & \checkmark & \checkmark    \\
\bottomrule
\end{tabular}
}
\end{table}

The HCLP can be seen as a hub location problem, since it aims to locate  some \emph{privileged} facilities, i.e., hubs, and to connect them together with served demand-nodes through a hub-and-spoke network. It differs, however, from most hub location problems in the literature in how demand is allocated and modelled, and in the topology induced by the hubs network.
In the HCLP, demand is captured if it is placed within a certain distance (coverage radius) of an activated hub. Thus, the HCLP can also be seen as a covering location problem \cite{GarciaMarin2nded}, whose main characteristic is that demand is \emph{captured} (covered) only if it is sufficiently close to open facilities. In our case, the covering distance, $\delta$, is exogenously defined, as in covering problems. 

One of the main differences between classical location covering problems and the HCLP is that, in the latter, the network induced by the activated facilities (hubs) must be connected. This further relates the HCLP to the location areas that impose the connectivity of the activated facilities as, for instance, hub location and  connected facility location.

Hub location is an active research area within Location Analysis as shown by the increasing number of works in the literature addressing these problems. Recent surveys can be found, for example, in \citep{Contreras2019, Alumur2021, OKelly2025}. Hubs are special facilities typically used to consolidate, sort, and distribute flows between multiple origin-destination (OD) pairs, often resulting in economies of scale.  Most HLPs focus on the design of hub-and-spoke networks with specific characteristics.

Early work in the area traces back to the seminal work of Goldman \citep{Goldman1969}, the problem statement and quadratic formulations by O'Kelly \citep{OKelly1986, OKelly1987}, and the linear reformulations developed by Campbell in analogy with median, center, covering, and fixed-charge facility location problems \citep{Campbell1992,Campbell1994}.

The HCLP is closely related to many HLPs and, in particular, to hub-covering problems, first introduced by Campbell \cite{Campbell1994}. As indicated in that paper, in the context of hub location, the notion of coverage may have several interpretations, depending on whether it is associated with OD pairs or with individual nodes. In the first case, an OD pair is covered if the total cost (time) of the path that connects its origin and destination through the backbone network does not exceed a specified value. Some HLPs dealing with this type of coverage have been studied in \cite{Wagner2008,Hwang2012,Tan2007,Kara2003,Alumur2009b}. In the latter case, a node is covered if the cost/time of the access/distribution arcs connecting the node to the backbone network does not exceed a specified value. Some works in which this type of coverage is considered are \cite{Lowe2013,Karimi2011,Peker2015,Sener2023}.

In this work, we combine both types of coverage. The first type is considered by imposing that, in the backbone network, the maximum travel time between any two activated hubs does not exceed a pre-specified time limit $T_{max}$. The second type of coverage is considered by imposing a maximum distance $\delta$ between any served demand node and the backbone network. Similar to covering facility location problems, with this type of coverage, it is no longer necessary to explicitly allocate covered nodes to activated hubs. To the best of our knowledge, both types of coverage have not yet been jointly considered in the literature.

Concerning the topology of the solution network, the HCLP does not require the backbone network to be complete, although connectivity is enforced. An increasing amount of research has focused on HLPs with incomplete backbone networks (see, e.g., \citep{Alumur2009, Davari2013, Calik2009, Sa2018}), which typically incorporate setup costs for activated interhub links. Even if, typically, optimal backbone networks of the HCLP tend to have very few arcs, it cannot be framed as a HLP with a specific-structure solution network such as a tree \citep{Contreras2010} or a line \citep{Sa2015}, given that cycles are allowed.

Similarities can also be found between the HCLP and the connected facility location problem (ConFL). The ConFL is a general network design problem that combines the Steiner tree and facility location problems. In particular, ConFL is the problem of determining a minimum-cost solution consisting of a set of facilities used to serve a given demand associated with a given set of users, and connecting the open facilities using a Steiner tree. The cost of a ConFL solution comprises the edge costs of building the Steiner tree, the cost of allocating users to facilities, and the cost of opening facilities.
The ConFL was first introduced in \citep{KargerMinkoff2000} as a probabilistic version of the rooted Steiner Tree Problem and has received considerable attention in the literature. Several variants and various solution approaches have been proposed (see, e.g., \citep{Ljubic2020}).  Applications of the ConFL arise in diverse areas such as telecommunications \citep{Aruselvan2019, Grotschel-et-al-2014, Leitner-et-al-2013, Ljubic-et-al-2006}, information/content distribution networks \citep{KargerMinkoff2000, KrickRacke2003}, and emergency management \citep{ Zhu-et-al-2012}. 

There are two main differences between the HCLP and the ConFLs. One is the maximum travel time between facilities, which is present in the HCLP but is not enforced in the ConFL. The other one is that the HCLP does not assume that the backbone network is a (Steiner) tree.

From the above literature review, we conclude that the HCLP is a novel model closely related to other models studied in Location Analysis, yet clearly differentiated from them.

\section{Problem definition}\label{sec:prob-def}

Consider a complete network $N=(V, A)$, with node set $V=\{1, 2,\dots, n\}$ and arc set $A=\{(i, j): i, j\in V\}$. 
Let $E=\{ij: (i, j)\in A\}$ denote the set of (undirected) edges underlying the arc set $A$. To alleviate notation, in the remainder of this paper, any edge $ij\in E$ with $i, j\in V$, $i<j$ will be indistinctively denoted by $ji$.
Also, for a given subset of nodes $S\subset V$, we denote by $E[S]=\{e=ij\in E: i,j\in S \}$ the set of edges with both endnodes in $S$. 
Associated with each arc $(i,j)\in A$, $d_{ij}$ denotes the distance from $i$ to $j$, and $t_{ij}$ the time required to traverse the arc in the transportation system. We assume that both $d_{ij}$ and $t_{ij}$ are symmetric and satisfy the triangle inequality.  

We assume, without loss of generality, that all nodes are potential hub locations in the transportation network. Each activated hub node $k\in V$ incurs a setup cost $G_k$. Furthermore, each activated hub $k\in V$ enables access to the transportation network to potential users $i\in V$ located at non-activated nodes within a given \emph{coverage radius} from $k$. The coverage radius is denoted by $\delta$ and is measured relative to the distances $d_{ik}$. For each potential user $i\in V$ we denote by $H_i$ the set of hubs that could enable service to this user, that is, $H_i=\{k\in V\,:\, d_{ik}\leq \delta\}$. Potential users $i\in V$ that are not within the coverage radius of any activated hub will be referred to as \textit{uncovered} (or \textit{non-covered}); such users incur a penalty $P_i>0$.

Edges connecting activated hub nodes can also be activated. Then their underlying arcs can be traversed in either direction. Each activated edge $ij\in E$ incurs a setup cost $F_{ij}$. The interhub network must be designed in such a way that the travel time between any two activated hubs does not exceed a given  {\it maximum travel time}, $T_{max}$.
In the following, pairs of potential hub nodes whose travel time exceed $T_{max}$ will be referred to as \emph{incompatible pairs}. Moreover, we will say that a subset of potential hubs $H \subseteq V$ is \textit{compatible} if any pair of them can be connected in a maximum travel time $T_{max}$, i.e.,  $t_{km}\leq T_{max}$ for all $k, m\in H$.

The HCLP is to determine, 
\begin{itemize}
\item A set $H \subset V$ of hub nodes to activate.
\item A set of hub edges to activate, $E_H\subseteq E[H]$, with underlying set of arcs $A_H$, 
such that, in the network $(H, A_H)$, the total travel time between any two activated hubs does not exceed $T_{max}$.
\item The set of potential users not covered by the solution network $(H, A_H)$: $U_H=\left\{i\in V\setminus H: \min_{k\in H}\{d_{ik}\}>\delta\right\},$
\end{itemize}

\noindent such that the sum of the following costs is minimized:
\begin{itemize}
\item Setup costs of activated hubs: $\sum_{k\in H}G_k$.
\item Setup cost of activated interhub edges: $\sum_{km\in E_H}F_{km}$.
\item Penalty of \textit{uncovered} users: $\sum_{i\in U_H}P_i$.
\end{itemize}

\begin{figure}[htb]
    \centering
    \includegraphics[width=0.4\linewidth]{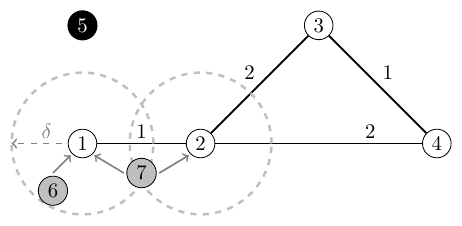}
    \caption{Illustrative example.}
    \label{fig:ejemplo}
\end{figure}

The solution network $(H, E_H)$ will also be referred to as the \textit{backbone network}. Figure~\ref{fig:ejemplo} illustrates a solution for an HCLP instance defined on a network consisting of seven nodes, with coverage radius $\delta$, and $T_{max}=3$. In this solution, the set of activated hubs is $H=\{1, 2, 3, 4\}$, and the set of activated interhub edges $E_H=\{12,\, 23,\, 24,\, 34\}$.
The travel times of the activated interhub connections are indicated next to the edges. For clarity, only the coverage areas of hubs 1 and 2 are shown as dotted circles, although hubs 3 and 4 also have their own coverage areas. As can be seen, nodes 6 and 7 are within a distance of at most $\delta$ of some activated hub, so users at those nodes can access the network: node 6 through hub 1, and node 7 through both hubs 1 and 2. In contrast, node 5 remains uncovered and will incur a penalty $P_5$. 
Note that the four edges in $E_H$ are needed to guarantee that the routing time between any pair of hubs does not exceed $T_{max}$. For instance, if interhub edge $34$ were not activated, the route from $3$ to $4$ would have to traverse the two arcs $(3,2)$ and $(2,4)$, implying a total travel time of $4>T_{max}$. Therefore, this example shows that to achieve travel-time feasibility, the backbone network can contain cycles. 

\subsection{Special cases of the HCLP}

To the best of our knowledge, the HCLP is a novel model, which has not been studied in the literature. However, it includes, as special cases, some well-known problems, as we discuss below.

\begin{itemize}
\item When $T_{max}=\infty$, the total travel constraint can be ignored. Then, \text{extreme} values for the interhub activation costs produce the following two special cases:

\begin{itemize}
\item $F_{ij}=0$ for all $ij\in E$. Then, similarly to many fundamental hub location models, there is an optimal solution where the backbone network $(H, E_H)$ is a complete graph. Hence, since the HCLP involves no routing costs, this special case involves neither network design decisions nor routing decisions, so it reduces to the problem of determining a set of facilities (hubs) to activate together with a set of \textit{captured} (covered) users within the coverage radius of the activated facilities. Such a problem belongs to the class of Covering Location Problems (see, e.g. \citet{GarciaMarin2nded}).

\item $F_{ij}=\infty$ for all $ij\in E$. Then, no interhub edge will be activated in any optimal solution. Thus, optimal solutions will take the form of a star topology, consisting of one single activated facility (hub) together with a set of \textit{captured} users within the coverage radius of the activated facility. Such a problem is a simple version of a Covering Location Problem mentioned above (see, again~\citet{GarciaMarin2nded}).
\end{itemize}

\item When $\delta=\infty$, then all users will be captured by any solution network. Thus, in this case the HCLP involves solely network design decisions (see., e.g., \cite{crainic2021network}), specifically, those focusing on determining a backbone network of minimum total setup cost (for the activated hubs and interhub edges) such that the traveling time between any two hubs does not exceed $T_{max}$.
\end{itemize}

In what follows, to avoid trivial solutions with only one hub activated, we assume that the backbone network contains at least one interhub link.

\section{Formulation with four-index variables ($F_1$)}\label{sec:4i}

In this section, we develop a first MILP formulation for the HCLP. In addition to binary variables that identify the activated hub nodes, the activated interhub links, and the users within the coverage radius of some activated hub, the formulation uses four-index variables to identify, for each pair of activated hubs, a path connecting them in the bacckbone network whose total travel time does not exceed the time limit $T_{max}$. 

The decision variables are the following: 
\begin{itemize}		
	\item $y_{k}$: binary variable that takes the value 1 if and only if a hub is opened at node $k \in V$, and 0 otherwise.
    \item $v_i$: binary variable that takes the value 1 if and only if user $i\in V$ is non-covered, and 0 otherwise.
\item $w_{km}$:  binary variable that takes the value 1 if and only if the pair of hubs $k,m\in V$ with $k< m$ are both activated, and 0 otherwise.
	\item $X_{km}$: binary variable that takes the value 1 if and only if the edge $km\in E$ is activated, and 0 otherwise.
	\item $x_{ij}^{km}$: binary variable that takes the value 1 if and only if the path connecting the pair of hubs $k,m\in V: k\neq m$ uses arc $(i, j)\in A$, and 0 otherwise.
\end{itemize}

The formulation is as follows.
\begin{subequations}
\begin{align}
\left(F_{1}\right) \quad  \text{min } \, & \sum_{k\in V}G_k y_k+\sum_{km\in E}F_{km}X_{km}  + \sum_{i\in V} P_i\,v_i&&  \label{FO}\\
\qquad \text{s.t.} \quad %
& \sum_{k\in H_i} y_k  + v_i\geq 1 \quad &&  i\in V   \label{C_service}\\
& w_{km}\leq y_k && k,m\in V :k< m  \label{C_wy_1}\\
& w_{km}\leq y_m && k,m\in V :k< m \label{C_wy_2}\\
& y_k+y_m\leq w_{km}+1&& k,m\in V :k< m  \label{C_wy_3}\\
& X_{km}\leq w_{km} && km\in E \label{C_Xw}\\
& \sum_{km\in E} X_{km}  \geq 1 \quad &&  \label{C_twohubs}\\
& \sum_{j\in V\setminus\{i\}}x^{km}_{ij}-\sum_{j\in V\setminus\{i\}}x^{km}_{ji}=\left\{
\begin{array}{ll}
w_{\{km\}} & k=i, \\
-w_{\{km\}} & m=i, \\
0 & \mbox{o/w,}
\end{array}
\right. &&k,m\in V :k\neq m, i\in V\label{C_flow}\\
& x_{ij}^{km}+x_{ji}^{km}\leq X_{ij} && k,m\in V :k\neq m , ij\in E\label{C_Xx}\\
& \sum_{(i,j)\in A}t_{ij}\,x^{km}_{ij}\leq T_{max}\,w_{km} \quad && k,m\in V :k\neq m  \label{C:Tmax}\\
& y_{k} \in \left\{0,1\right\} \quad && k \in V  \label{Dom:y}\\
& v_i\in \left\{0,1\right\} \quad && i\in V \label{Dom:s}\\
& w_{km} \in \left\{0,1\right\} \quad && k,m\in V :k< m   \label{Dom:w} \\
& X_{km} \in \left\{0,1\right\} \quad && km\in E  \label{Dom:X} \\ 
& x_{ij}^{km} \in \left\{0,1\right\} \quad && k,m\in V :k\neq m , (i,j)\in A.  \label{Dom:x}
\end{align}
\end{subequations}

The objective function minimizes the sum of the activation cost of the network (hubs and interhub edges) plus a penalty term for uncovered nodes. Constraints \eqref{C_service} identify as non-covered users all potential users that are not in the coverage radius of at least one activated hub. Because of their positive coefficients in the objective function, in an optimal solution $v_i=0$ for all covered users. The linearization of the products $w_{km}=y_k\,y_m$ is defined in \eqref{C_wy_1}-\eqref{C_wy_3} for all pairs of potential hubs.
The relationship between the variables $X$  and $w$  is established in \eqref{C_Xw}. These constraints model the logic that no interhub edge can be activated unless its two endnodes are activated hubs (i.e. $X_{km}\leq y_k$ and $X_{km}\leq y_k$). However,  because of \eqref{C_wy_1}-\eqref{C_wy_2},  they are tighter than the \textit{direct} modeling of this logic, so we use them in the formulation. The purpose of constraint \eqref{C_twohubs} is to ensure that the backbone network contains at least one interhub edge. 
Constraints \eqref{C_flow} determine the arcs of the paths connecting pairs of activated hubs. By constraints \eqref{C_Xx}, such paths may only use activated interhub links, and by  \eqref{C:Tmax} their total travel time does not exceed the maximum allowed duration $T_{max}$. Finally, \eqref{Dom:y}-\eqref{Dom:x} define the domain of the decision variables.

The feasible solution given in Figure~\ref{fig:ejemplo} corresponds to the following non-zero variables for $F_1$. Regarding activated hubs and non-covered users, $y_1=y_2=y_3=y_4=1$ and $v_5=1$ respectively. The interhub connections are associated with variables $X_{12}=X_{23}=X_{24}=X_{34}=1$. The activated pairs of hubs and their corresponding paths that do not exceed the maximum travel time, are  described by the following $w$ and $x$ variables, respectively: 
$w_{12}=1$, $x^{12}_{12} =1$; $w_{13}=1$, $x^{13}_{12}= x^{13}_{23}=1$; $w_{14}=1$, $x^{14}_{12} = x^{14}_{24}=1$; $w_{23}=1$, $x^{23}_{23} =1$; $w_{24}=1$, $x^{24}_{24} =1$; and, $w_{34}=1$,  $x^{34}_{34} =1$.

\subsection{Feasibility and optimality conditions for $F_1$}
Several feasibility and optimality conditions apply to formulation $F_1$. Below we mention those that have proven to have some effect in preliminary computational testing. 

\begin{enumerate}
\item[F1] A necessary condition to have a feasible solution is that the set of activated hubs be compatible.  
Therefore, for all $k,m\in V$ such that \,$t_{km}>T_{max}$, $w_{km}=0$, $X_{km}=0$, $x_{ij}^{km}=0$ for all $(i,j)\in A$, and $x_{km}^{st}=0$ for all $s,t\in H, s\neq t$.

\item[O1] As each path begins and ends at the designated hubs, there is an optimal solution where $x_{ik}^{km}=x_{mj}^{km}=0$ for all $k,m\in V:k\neq m$ and $(i,k), (m,j)\in A$. 

\item[O2] Since we assume a complete network and symmetric travel times, there is an optimal solution where $x_{ij}^{km}=x_{ji}^{mk}$ for all $k,m\in V: k<m,\, (i,j)\in A$. Therefore, the number of path variables can be reduced to $x_{ij}^{km}$ for $k,m\in V: k<m,\, (i,j)\in A$.
\end{enumerate}

\subsection{Valid inequalities}

We have studied a variety of valid inequalities for formulation $F_1$. Below we include those that, in preliminary testing,  had some effect on the linear relaxation (LP)  bound or the computing times:

\begin{enumerate}
\item The backbone network has to be connected; otherwise, paths connecting pairs of hubs in different components would not exist. Hence, the number of edges of the backbone network has to be greater than or equal to the number of activate hubs minus one. 
\begin{align} 
& \sum_{km\in E} X_{km} \geq \sum_{k\in V} y_k -1. && \label{eq:connection}
\end{align}
Note that the equality cannot be imposed because the backbone network may contain cycles, to allow feasible paths (regarding time limits) connecting pairs of activated hubs. Note also that, when $\sum_{k\in V} y_k>2$, i.e. more than two hubs are activated, these inequalities dominate constraints \eqref{C_twohubs}.

\item Incompatible pairs of nodes cannot be activated as hubs simultaneously. Therefore, in addition to the variable fixing indicated above, the following inequalities can be used to reinforce the LP bound:
\begin{align}
& y_k + y_m\leq 1 &&k,m\in V: t_{km}>T_{max}. \label{vi_clique_1}   
\end{align}

Inequalities \eqref{vi_clique_1} can be generalized to pairwise-incompatible subsets of potential hubs. In particular, let $Q\subset V$ such that $t_{km}>T_{max}$ for all $k,m\in Q$. Then, the following \textit{clique} inequality is valid
\begin{align}
& \sum_{k\in Q}y_k\leq 1. && \label{vi_clique_Gen}   
\end{align}

\end{enumerate}

\section{Formulation with three-index variables ($F_2$)}\label{sec:3i}

In this section we develop a three-index MILP formulation for the HCLP, which will be referred to as $F_2$. The strategic decisions are modeled using the same variables as $F_1$, so the variables $X$, $y$, $v$, and $w$ are the same as above, but $F_2$ uses a different approach to connect activated hubs and for the time-limit constraints.
In particular, $F_2$ builds a tree rooted in each activated hub, which \textit{spans} all other activated hubs. Now, the backbone network is determined by all the edges underlying the arcs of the spanning trees.
Since any tree contains one single path connecting its root and any of its nodes, the feasibility of the backbone network can be established using some MTZ-type constraints \cite{MTZ} that compute the travel times from the root to any other activated hub, and imposing that these times do not exceed $T_{max}$.
$F_2$ replaces the four-index path variables $x^{km}_{ij}$ of $F_1$ by two sets of variables, one binary and one continuous. The binary variables have three indices and model the (directed) spanning trees rooted at the activated hubs, whereas the continuous variables have two indices and determine the travel times from each tree root to each activated hub.
They are defined as follows:
\begin{itemize}		
	\item $x^s_{km}$: binary variable that takes the value 1 if and only if arc $(k,m)\in A$ is in the tree rooted at hub $s\in V$, and 0 otherwise. 
	\item $u_{sk}$: continuous variable measuring the travel time from origin hub $s\in V$ to destination hub $k\in V$ through the tree rooted in $s$. Note that, by the symmetry of travel times through the arcs, $u_{sk}$ also determines the travel time from $k$ to $s$. 
\end{itemize}

The formulation $F_2$ is detailed below. 
\begin{subequations}
\allowdisplaybreaks 
\begin{align}
(F_{2}) \quad  \text{min } \quad & \sum_{k\in V}G_k y_k+ 
\sum_{km\in E}F_{km}X_{km}  + \sum_{i\in V} P_i\,v_i&&  \label{FO2}\\
\qquad \text{s.t.} \quad & \sum_{k\in H_i} y_k  + v_i\geq 1 \quad &&  i\in V   \label{C2_service}\\
& w_{km}\leq y_k && k,m\in V: k< m   \label{C2_wy_1}\\
& w_{km}\leq y_m && k,m\in V: k< m   \label{C2_wy_2}\\
& y_k+y_m\leq w_{km} +1&& k,m\in V: k< m    \label{C2_wy_3}\\
& X_{km}\leq w_{km} && km\in  E \label{C2_Xw}\\ 
& \sum_{km\in E} X_{km}  \geq 1 \quad &&  \label{C2_twohubs}\\ 
& x^s_{km} + x^s_{mk}\leq X_{km} &&  s\in V,\, km\in E\label{C2_Xx}\\ 
& \sum_{(s,m)\in A} x^s_{sm} \geq y_s && s\in V \label{C2_tree-root}\\
& \sum_{(k,m)\in A} x^s_{km} = w_{sm}  && s,m\in V: s\neq m \label{C2_xw}\\
&  u_{ss} = 0 \quad && s\in V  \label{C2_uzero}\\
&  u_{sk}\geq u_{sm} + t_{mk} - (T_{max}+ t_{mk})\, (1-x^s_{mk}) \quad && s,k\in V: s\neq k,\, (m,k)\in A \label{C2_MTZ}\\ 
&  u_{sk}\leq t_{sk} + (T_{max} - t_{sk}) \,(1-X_{sk}) \quad && s,k\in V: s\neq k  \label{C2_udirect}\\
& y_{k} \in \left\{0,1\right\} \quad && k \in V  \label{Dom2:y}\\
& v_i\in \left\{0,1\right\} \quad && i\in V \label{Dom2:s}\\
& w_{km} \in \left\{0,1\right\} \quad && k,m\in V: k< m  \label{Dom2:w} \\
& X_{km} \in \left\{0,1\right\} \quad && km\in E  \label{Dom2:X} \\ 
& x_{km}^{s} \in \left\{0,1\right\} \quad && s\in V,\,(k,m)\in A\label{Dom2:x}\\ 
& u_{sk} \geq 0\quad && s,k\in V.\label{Dom2:u}
\end{align}
\end{subequations}

The objective function \eqref{FO2} and constraints \eqref{C2_service}--\eqref{C2_twohubs} are, respectively, the same as  \eqref{FO} and \eqref{C_service}-\eqref{C_twohubs} in formulation $F_1$. Constraints \eqref{C2_Xx} guarantee that no arc is used in a tree unless the underlying interhub edge has been activated. Constraints \eqref{C2_tree-root} and \eqref{C2_xw} define the tree structure associated with each activated hub: while constraints \eqref{C2_tree-root} impose that at least one arc leaves  node $s$, if  it is is activated as a hub, by constraints \eqref{C2_xw} exactly one arc of the tree rooted in hub $s$ enters any other activated hub. Constraints \eqref{C2_uzero} and \eqref{C2_MTZ} determine the values of variables  $u_{sk}$: while constraints \eqref{C2_uzero} set the initial travel times, constraints \eqref{C2_MTZ} regulate their values taking into account the arcs used in the rooted trees, using the same rationale as in classical MTZ-type of constraints \cite{MTZ}. To guarantee the feasibility of solutions, in \eqref{C2_udirect}  the upper bounds for the values of these variables  are set to $T_{max}$.

The feasible solution given in Figure~\ref{fig:ejemplo} is associated with the following values of the variables of $F_2$. Variables $X$, $y$, $w$ and $v$ take the same values as in $F_1$. The non-zero variables that represent the interhub arcs of the trees rooted at the activated hubs ($x^s_{km}$), which allow to determine the travel times from the roots to other activated hubs ($u_{sk}$) are detailed in Table~\ref{table:example_F2-2}. Note that the values of the $u$-variables involving hubs located at leaves of the trees always correspond to travel times of feasible paths that do not exceed $T_{max}$, although they do not necessarily correspond to minimum travel times in the backbone network.

\begin{small}
\begin{table}[htbp]
\begin{center}
\caption{Values of $F_2$ variables $x^s_{km}$ and $u_{sk}$ for the  example of Figure~\ref{fig:ejemplo}.}
\label{table:example_F2-2}
\begin{tabular}{cll}
\toprule
Hub & Tree variables & Time variables\\ 
\midrule
1 & $x^1_{12}=x^1_{23}=x^1_{24}=1$ & $u_{12}=1$, $u_{13}=3$, $u_{14}=3$\\
2 & $x^2_{21}=x^2_{23}=x^2_{24}=1$ & $u_{21}=1$, $u_{23}=2$, $u_{24}=2$\\
3 & $x^3_{32}=x^3_{21}=x^3_{34}=1$ & $u_{32}=2$, $u_{31}=3$, $u_{34}=1$\\
4 & $x^4_{42}=x^4_{21}=x^4_{43}=1$ & $u_{42}=2$, $u_{41}=3$, $u_{43}=1$\\
\bottomrule
\end{tabular}
\end{center}
\end{table}
\end{small}
\vspace{-0.50cm}

\subsection{Feasibility and optimality conditions for $F_2$}

Feasibility and optimality conditions have also been studied for formulation $F_2$. Below we mention those that have proven to have some effect in preliminary computational testing. 

\begin{enumerate}
\item[F2]  Similarly to $F_1$, no pair of nodes  $k,m\in V$ such that \,$t_{km}>T_{max}$ can be activated simultaneously in any feasible solution. Hence, as before, $w_{km}=0$ and $X_{km}=0$ for all $k,m\in V$ such that \,$t_{km}>T_{max}$. Moreover, $x_{km}^{s}=0$ for all $s\in V$ and $u_{km}=0$.

\item[O3] There is an optimal solution where $x^s_{ks}=0$ for all $s\in V,\, (k,s)\in A$. 

\item[O4] As we assume that the input network is complete, and symmetry in distances and travel time, the number of variables $u_{sk}$ can be reduced as follows. For any $s,k\in V$, any feasible solution will satisfy that $u_{sk} = u_{ks}$, therefore we will only consider variables $u_{sk}$ for $s,k\in V\,:\,s<k$.

\end{enumerate}

\subsection{Valid inequalities}
Inequalities \eqref{eq:connection}, \eqref{vi_clique_1} and \eqref{vi_clique_Gen} are also valid for $F_2$. After some preliminary testing, we have observed they also have an effect on its LP bound. 

\section{Two-phase solution algorithm}\label{sec:MathH}

In this section, we present a two-phase solution algorithm for the HCLP. The rationale for this matheuristic is to decompose the formulation $F_2$ into two independent subproblems: LC, which focuses only on location and covering decisions, and $ND(H)$, which yields optimal network design decisions for a given set of activated hubs, $H\subseteq V$. Indeed, by combining the outcome of both phases, in each iteration, we obtain a feasible HCLP solution. Both subproblems are detailed below, along with the proposed matheuristic. 

\subsection{Phase 1: Location and Covering subproblem}
The \textit{Location and Covering} subproblem, LC, determines a set of compatible hubs and a set of non-covered users that minimize the sum of hub activation costs and the penalties for non-served users. Its objective function value will be denoted by $z_{LC}$. Subproblem LC includes a set of constraints \textit{inherited} from $F_2$, which involve variables $y$ and $v$. Since LC is solved once per iteration of the matheuristic, this subproblem also includes additional sets of constraints that are activated only after the first iteration and prevent the same set of hubs from being opened across iterations. These sets of constraints will be dynamically updated as the overall matheuristic evolves and guide the search towards new sets of activated hubs. 
The formulation is as follows:
\begin{subequations}
\allowdisplaybreaks 
\begin{align}
\text{LC}(\rho, \mathcal P, H^*, z_{\text{HCLP}}^*) \quad 
\text{min }  \quad & z_{LC}=\sum_{k\in V}G_k y_k + \sum_{i \in V} P_i v_i&&  \label{LCP_FO}\\
\qquad\text{s.t.}\quad 
& \sum_{k\in H_i} y_k + v_i \geq 1 \quad &&  i\in V   \label{LCP_cover}\\
& \sum_{k\in V} y_k  \geq 2 \quad &&  \label{LCP_twohubs}\\
& y_k + y_m \leq 1 &&  k,m \in V : t_{km} > T_{max} \label{LCP_Cliques_1}\\
& \displaystyle \sum_{k \in V\backslash P }y_{k} + \displaystyle \sum_{k \in P} (1-y_k) \geq 1  && P \in \mathcal P \label{LCP_Pool}\\
& \displaystyle \sum_{k \in V\backslash H^*} y_{k} + \displaystyle \sum_{k \in H^*} (1-y_k) = \rho && \label{LCP_radius}\\
& \sum_{k\in V}G_k y_k + \sum_{i \in V} P_i v_i< z_{\text{HCLP}}^* - \widehat{F}  && \label{LCP_zbest}\\
& y_{k} \in \left\{0,1\right\} \quad && k \in V  \label{LCP_Dom:y}\\
& v_i \in \left\{0,1\right\} \quad && i \in V. \label{LCP_Dom:v}
\end{align}
\end{subequations}

The objective function~\eqref{LCP_FO} involves the setup costs for the activated hubs plus the penalties for non-covered nodes. Constraints~\eqref{LCP_cover}  are the inequalities \eqref{C2_service} of  $F_2$,  which identify the non-covered nodes. Inequalities~\eqref{LCP_twohubs} express the logic of constraints \eqref{C2_twohubs} in terms of the $y$-variables, whereas \eqref{LCP_Cliques_1} are the inequalities \eqref{vi_clique_1}, which discard incompatible pairs of hubs.

Constraints~\eqref{LCP_Pool}-\eqref{LCP_zbest}, which depend on the set of parameters  $(\mathcal P, z_{\text{HCLP}}^*, H^*, \rho)$,  guide the process after the first iteration. The parameters are updated through an iterative procedure as explained below. The set $\mathcal P$ denotes a pool consisting of the sets of hubs activated in previous iterations (initially, $\mathcal P=\emptyset$). $z_{\text{HCLP}}^*$ is the $F_2$ objective function value of the best-known HCLP  solution and $H^*$ its associated set of activated hubs.
Due to constraints \eqref{LCP_radius}, subproblem $LC$ can be interpreted as a local search in a neighborhood of the currently best-known HCLP solution,  where the parameter $\rho$ is the radius (\textit{depth}) of the neighborhood. 
The purpose of this local search is to obtain \textit{promising} sets of hubs in the current neighborhood. Thus, solving subproblem $LC$ for a given set $H^*$ and varying values of $\rho$ can be viewed as a variable neighborhood search, in which the subproblems associated with small values of $\rho$  are \textit{intensification} moves, whereas those associated with large values of $\rho$ are \textit{diversification} ones.

Overall, constraints~\eqref{LCP_Pool}-\eqref{LCP_radius} are used to prevent repetitions and to balance 
\textit{intensification} and \textit{diversification} relative to best-known solutions. 
Specifically, constraints~\eqref{LCP_Pool} are no-good cuts that exclude sets of hubs that have already been explored. 
As explained, constraints~\eqref{LCP_radius} restrict the search space for potential sets of activated hubs to those within a neighborhood of radius $\rho$ of $H^*$, and produce solutions with exactly $\rho$ potential hubs whose \textit{status} is not the same as in $H^*$. Finally, constraints~\eqref{LCP_zbest} restrict the search to solutions with an objective function value for $z_{LC}$ strictly smaller than $z_{\text{HCLP}}^* - \widehat{F}$, where $\widehat{F}=\min\left\{F_{km}: km\in E \text{ s.t. $k$ and $m$ compatible}\right\}$.  Note that  $z_{\text{HCLP}}^* - \widehat{F}$ is an upper bound of the LC objective function value for any HCLP solution better than the current best-known.

We observe that the LC subproblem solved in the first iteration will be feasible (provided that the original HCLP is feasible), whereas constraints~\eqref{LCP_Pool}-\eqref{LCP_radius} may render the subproblem infeasible in subsequent iterations. For instance, when the current best-known solution is optimal, or when the explored neighborhood contains no solution better than the current incumbent. In such cases, the outcome of the algorithm will be an empty subset of activated hubs, $H=\emptyset$,  associated with an objective function value $z_{LC}=\sum_{i\in V}P_i$.

\subsection{Phase 2: Network Design subproblem}
Given a compatible set of hubs $H\subseteq V$, the Network Design subproblem $\text{ND}(H)$ aims to determine a feasible backbone network over $H$ of minimum value with respect to edge activation costs. We will denote by $z_{ND}$ its objective function value. The formulation is as follows:

\begin{subequations}
\begin{align}
  ND(H) \quad   
  \text{min } \quad & z_{ND}=\sum_{km\in E_{H}} F_{km}X_{km} &&  \label{ND_FO}\\
\qquad \text{s.t.} \quad 
& \sum_{km\in E_{H}} X_{km} \geq |H| -1 && \label{C2_connection_ND}\\
& \sum_{(k,m)\in A_{H}} x^s_{km} \leq |H|-1  && s\in H \label{C2_tree-arcs_ND}\\
& x^s_{km} + x^s_{mk}\leq X_{km} && s\in H;\,  km\in E_{H}\label{C2_Xx_ND}\\ 
& \sum_{(s,m)\in A_{H}} x^s_{sm} \geq 1  && s\in H \label{C2_tree-root_ND}\\
& \sum_{(k,m)\in A_{H}} x^s_{km} = 1  && s,m\in H: s\neq m \label{C2_xw_ND}\\
&  u_{ss} = 0 \quad && s\in H  \label{C2_uzero_ND}\\
&  u_{sk}\geq u_{sm} + t_{mk} - \left(T_{max}+ t_{mk}\right)\, (1-x^s_{mk}) \quad &&  s,k\in H: s\neq k,\, (m,k)\in A_{H}\label{C2_MTZ_ND}\\
&  u_{sk}\leq t_{sk} + T_{max}\left(1-X_{sk}\right) \quad && s,k\in H: s\neq k  \label{C2_udirect_ND}\\
& X_{km} \in \left\{0,1\right\} \quad && km\in E_{H}  \label{Dom2:X_ND} \\
& x_{km}^{s} \in \left\{0,1\right\} \quad && s\in H,\,(k,m)\in A_{H}\label{Dom2:x_ND}\\
& u_{sk} \geq 0\quad && s,k\in H.\label{Dom2:u_ND}
\end{align}
\end{subequations}  

The objective function~\eqref{ND_FO} only considers the setup costs of the interhub connections. 
Constraint~\eqref{C2_connection_ND} is the adaptation of the valid inequality \eqref{eq:connection}.
Constraints~\eqref{C2_tree-arcs_ND} establish the number of arcs in each tree. Now, since the hubs are already fixed, the number of arcs in each rooted tree can be fixed exactly to the number of hubs minus one. 
Constraints \eqref{C2_Xx_ND}--\eqref{C2_udirect_ND} are the adaptation of Constraints \eqref{C2_Xx}--\eqref{C2_udirect} from $F_2$ when restricted to the active set of hubs $H$ and interhub connections $E_H$. 

Indeed, the optimality conditions of $F_2$ also apply to ND. 
In particular, the symmetry can be applied to reduce the number of $u_{sk}$ variables considering only $s,k\in H\,:\, s<k$; moreover, $x_{ks}^s$ are also excluded from the formulation. Therefore, the number of variables and constraints is reduced.

\subsection{Matheuristic}
\begin{figure}[htbp]
    \centering
    \includegraphics[width=0.6\linewidth]{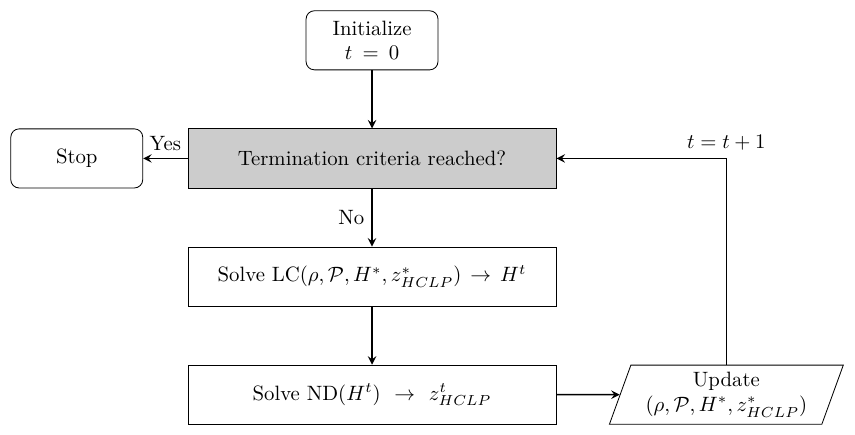}
    \caption{Matheuristic flowchart.}
    \label{fig:flowchart}
\end{figure}

Our matheuristic solves subproblems $\text{LC}(\rho, \mathcal P, H^*, z_{\text{HCLP}}^*)$ and $\text{ND}(H)$ in an iterative manner while applying a variable neighborhood search to obtain near-optimal solutions of the HCLP. A flowchart of the overall procedure is given in Figure~\ref{fig:flowchart}.  

At each iteration $t$, a compatible set of hubs $H^{t}$ is obtained by solving subproblem $\text{LC}(\rho, \mathcal P, H^*, z_{\text{HCLP}}^*)$.  If LC is infeasible, then $H^{t}=\emptyset$, $z_{LC}=\sum_{i\in V}P_i$,
and the neighborhood radius will be increased by one.

After solving subproblem LC, the set of activated hubs $H^{t}$ is added to the pool of solutions $\mathcal{P}$. 
Then, an optimal backbone network for the set of hubs $H^{t}$ is obtained by solving subproblem $ND(H^t)$, and the value $z_{\text{HCLP}}^{t}=z^t_{LC}+z^t_{ND}$ is computed. In case $z_{\text{HCLP}}^{t}$ is smaller than the value of the best-known solution, the incumbent is updated. Otherwise,  the neighborhood radius will be increased to enhance diversification. The procedure terminates when the time limit or the maximum neighborhood radius are reached.

The detailed pseudo-code is given in Algorithm~\ref{alg:matheur_pseudocode}. 
The input parameters are the maximum values for: computing time ($\tau_{max}$),  neighborhood radius ($\rho_{max}$), and number of solutions explored in the same neighborhood, i.e., number of iterations without improvement before increasing the depth of the neighborhood  
($\eta_{max}$). The algorithm's output is the objective function value of the best HCLP solution obtained within the time limit, $z_{\text{HCLP}}^*$. 

The algorithm starts in line 2, setting to 0 the initial values of the computing time ($\tau$) and algorithm counters: number of iterations ($t$) and number of iterations exploring the same neighborhood ($\eta$). Since the constraints that use the LC parameters, $(\mathcal P, z_{\text{HCLP}}^*, H^*, \rho, )$ are ignored in the first iteration, in lines 3-5  these parameters are initialized to \textit{fictitious} values that will be updated in the following iterations with the information corresponding to the HCLP solution obtained. In particular, we initially set $\mathcal{P}=\emptyset$, $z_{\text{HCLP}}^*=\sum_{i\in V} P_i$, $H^*=\emptyset$, and  $\rho=1$.
 
Once initialized, in lines 6-23 the algorithm continues as long as the termination criteria are not reached. In each iteration, subproblems LC and ND are solved in lines 8 and 9, respectively, and the parameters are updated accordingly: $\mathcal{P}$ in line 10, and $z_{\text{HCLP}}^{t}$ in line 11. When $z_{\text{HCLP}}^{t}$ improves the value of the best-kown solution, then $z_{\text{HCLP}}^{*}$ and $H^*$ are updated in lines 14 and 15, respectively, and the counters $\rho$ and $\eta$ are initialized. Otherwise, the counter for the number of iterations in the same neighborhood, increases in line 19, unless it has reached its limit. In such a case, it is initialized and the neighborhood radius is updated (lines 21-22).  The algorithm terminates either when the time limit has been reached or when all the neighborhoods with radius $1\leq \rho\leq\rho_{max}$ have been explored without improving the current incumbent solution.

\begin{small}
\SetKwComment{Comment}{$\triangleright$}{ }
\begin{algorithm}
\caption{Pseudo-code for the two-phase matheuristic method.}\label{alg:matheur_pseudocode}
\DontPrintSemicolon
\KwIn{$\tau_{max}$,\, $\rho_{max}$, $\eta_{max}$}
\KwOut{$z_{\text{HCLP}}^*$}
\SetKwBlock{Begin}{Procedure}{}
\Begin()
{
  $\tau, t, \eta \gets 0$\;
  $\rho \gets 1$\;
  $\mathcal P, H^* \gets \emptyset$\;
  $z_{\text{HCLP}}^* \gets \sum_{i\in V} P_i$\;
   \While{termination criteria not reached}
  {
        $t \gets t + 1$\;   
        $(H^{t}, z^{t}_{LC}) \gets \text{LC}(\rho, \mathcal{P}, H^*, z_{\text{HCLP}}^*)$
        \Comment*[r]{\text{Solve LC}}
         $z_{ND}^{t} \hspace{0.75cm}\gets \text{ND}(H^{t})$
            \Comment*[r]{\text{Solve ND}}                
            $z_{\text{HCLP}}^{t} \gets z^{t}_{LC} + z_{ND}^{t}$\;     
            $\mathcal P \gets \mathcal P \cup H^{t}$ \;
            \If{$\left(z_{\text{HCLP}}^{t} < z_{\text{HCLP}}^*\right)$}
            {                
                $z_{\text{HCLP}}^* \gets  z_{\text{HCLP}}^{t}$\Comment*[r]{\text{Update incumbent}}
                $H^* \gets H^{t}$\;
                $\eta \gets 0$\;
                $\rho \gets 1$\;
           }                 
            \ElseIf{$\left(\eta\leq\eta_{max}\right)$}
            {$\eta\gets \eta+1$ \Comment*[r]{\text{Intensify}}}
            \Else
            {   
                $\eta \gets 0$\;  
                {$\rho \gets \rho+1$\Comment*[r]{\text{Diversify}}}
            }                
            {\bf update} $\tau$\;
        }
  }
\end{algorithm}
\end{small}

\section{Computational experiments}\label{sec:compu}
In this section, we detail the computational experiments we have run and derive some managerial insights from the obtained solutions. The formulations and the algorithm described in the previous sections have been implemented using C++ and CPLEX 22.1, on CAI-URANIA\footnote{\href{https://supercomputacion.uca.es/cluster-de-supercomputacion/hardware-2/}{https://supercomputacion.uca.es/cluster-de-supercomputacion/hardware-2/}} computing cluster (2.6 GHz Intel Xeon-Platinum 8358 CPU). 
The branch-and-cut algorithm of CPLEX was restricted to run on a single thread with 100 GB of memory, to enable a better analysis of the results by excluding the parallel speedup.

This section is structured in the following parts. Firstly, we detail the benchmark instances generator. Then, we compare the performance of formulations $F_1$ and $F_2$, detail the performance of the matheuristic, and analyze the results of formulation $F_2$  using the matheuristic solution as a warm start. 
After that, we derive some managerial insights and perform a sensitivity analysis.
Finally, some public transport examples are given which emulate the use of real bus lines.

\subsection{Benchmark instances generator}\label{sec:instances-generator}
Preliminary testing has been carried out to produce instances with optimal solutions offering a suitable trade-off among the different elements. Instances\footnote{\href{https://github.com/carmenanadb/instances/tree/main/2026_HCLP}{https://github.com/carmenanadb/instances/tree/main/2026\_HCLP}} are derived from the well-known CAB data set\footnote{\href{https://www.researchgate.net/publication/269396247_cab100_mok}{https://www.researchgate.net/publication/269396247\_cab100\_mok}}, generating or adapting additional data when needed. From the $100$-node CAB instance, we use the distances between each pair of nodes $d_{ij}$, and we assume that travel times are proportional to distances, that is, $t_{ij} = \tau \, d_{ij}$ with $\tau = 0.1$ the conversion factor. Therefore, distances and travel times matrices are symmetrical, with zero diagonal, and their values satisfy the triangular inequality. 

Then, the setup costs are generated as follows. The hub activation costs are integer values rounded from those generated in \cite{Camargo2008} (also used in \cite{Alibeyg2016,Alibeyg2018,DominguezBravo2024}). The costs for activating an interhub edge are set as $F_{km} = \sigma\,d_{km}$ for each $km\in E_H$ with $\sigma = 100$. Finally, the penalty values for uncovered nodes are set equal to the activation costs of the nodes, i.e., $P_i = G_i$ for all $i\in V$. 

The parameters for maximum time and coverage radius are set as $T_{max}=230$ and $\delta=340$. These values satisfy that at least one node in the network is in the coverage radius of any other node, and that the entire network can not be covered by one single node. 

Instances with $|V|\in\{10,20,30,40,50,60,70,80,90,100\}$ are generated from the original $100$-node CAB network by selecting the $n$ nodes with more demand. Therefore, the nodes of the smaller instances will always be contained in the larger networks. In the following, we distinguish between {\it small instances} (up to $50$ nodes) and {\it large instances} (more than $50$ nodes). 

\begin{figure}[ht]
\begin{subfigure}[b]{0.35\textwidth}
\frame{\includegraphics[width=\textwidth]{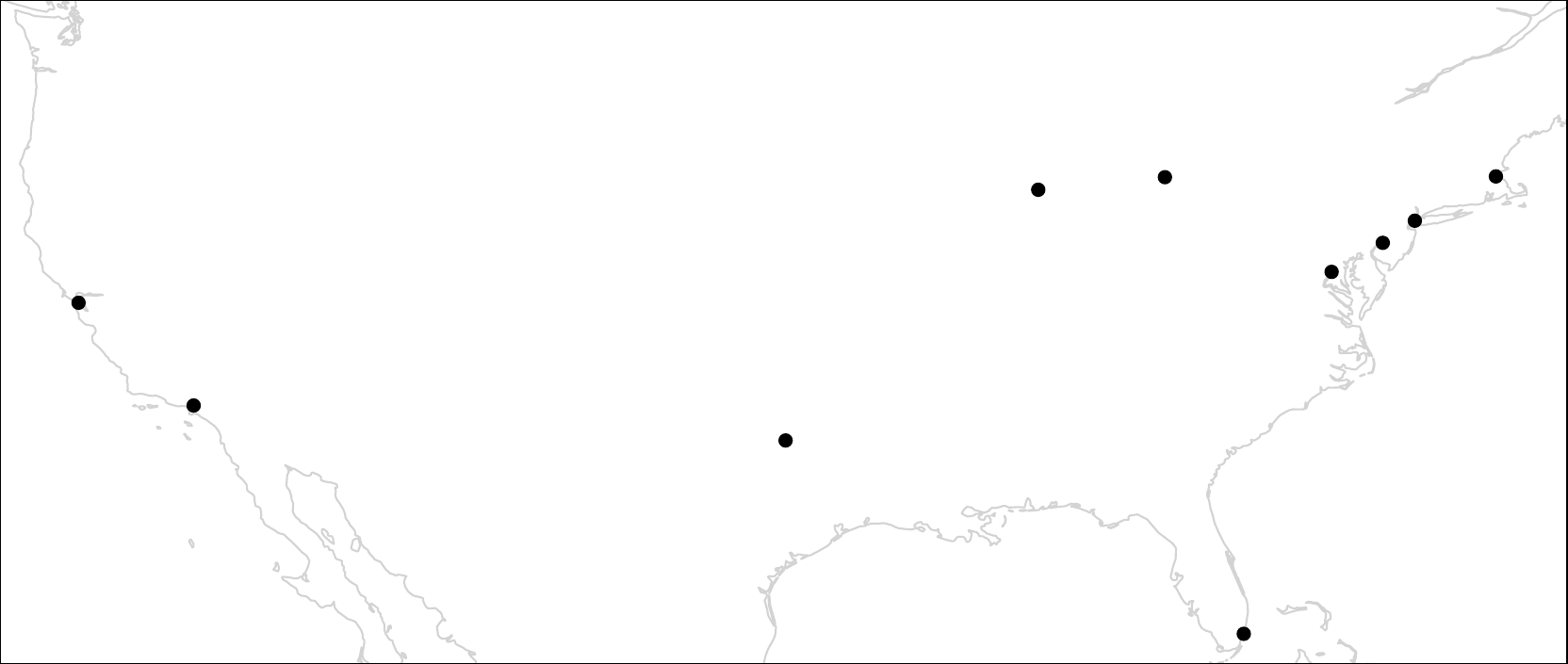}}
\caption{$|V|=10$}
\end{subfigure}
\begin{subfigure}[b]{0.35\textwidth}
\frame{\includegraphics[width=\textwidth]{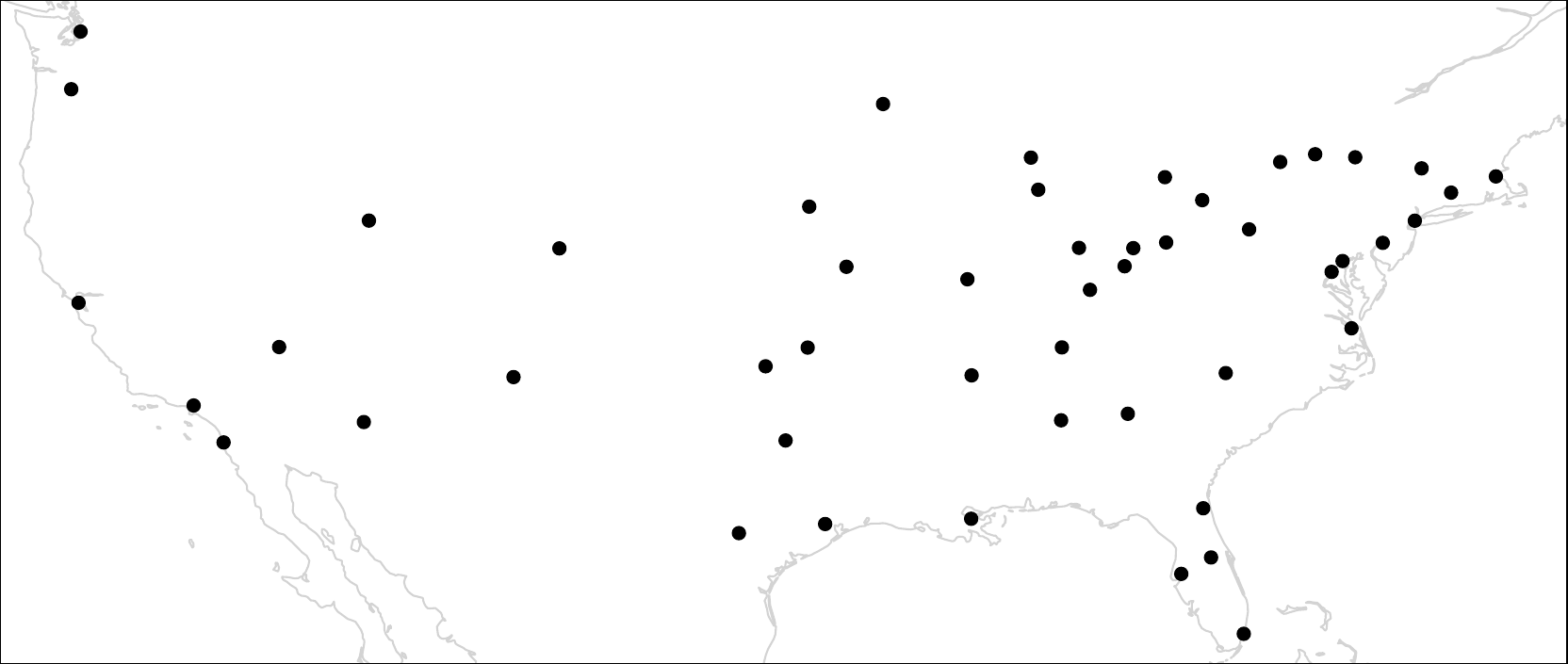}}
\caption{$|V|=50$}
\end{subfigure}
\begin{subfigure}[b]{0.35\textwidth}
\frame{\includegraphics[width=\textwidth]{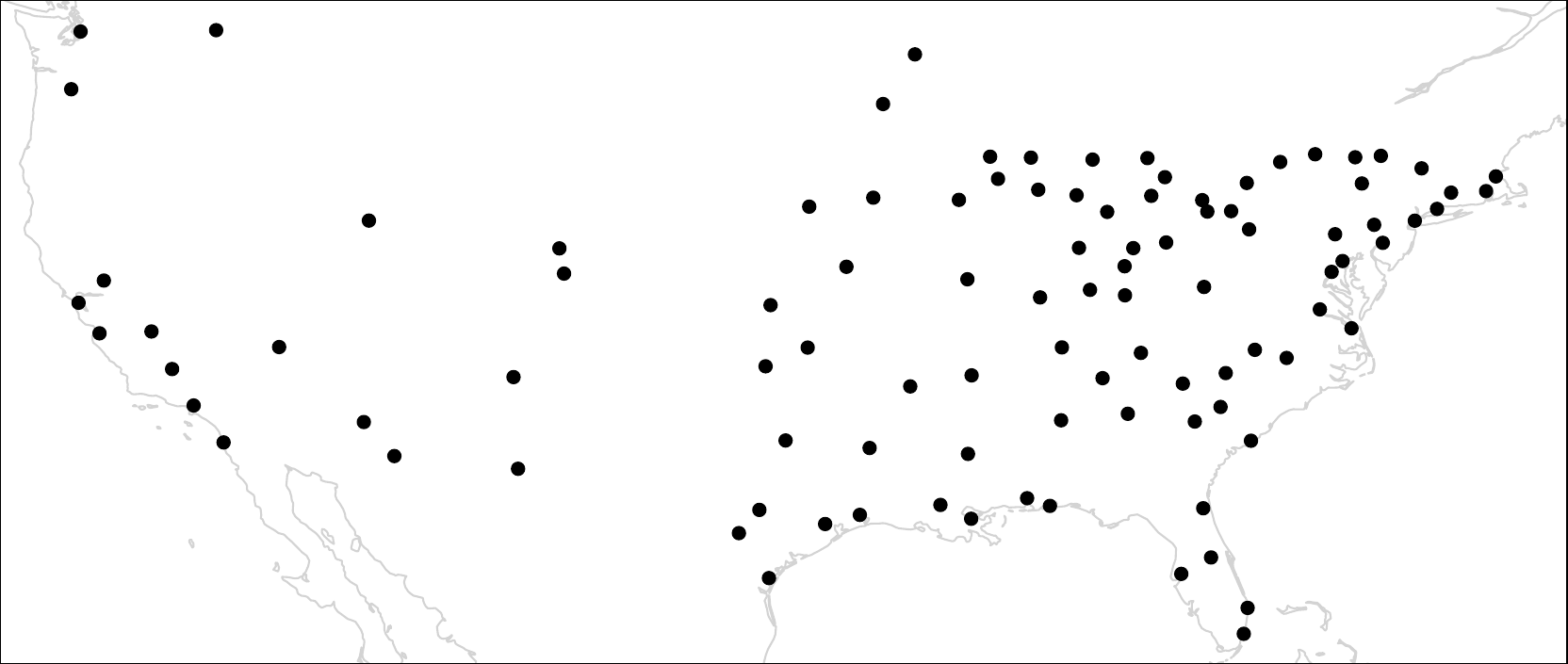}}
\caption{$|V|=100$}
\end{subfigure}
\caption{CAB network's nodes spatial distribution.}
\label{fig:CAB:nodes}
\end{figure}

Figure~\ref{fig:CAB:nodes} gives a general idea of the node spatial distribution for the different instances. Nodes located in both the west and east coast are already present in the smaller instance with 10 nodes. As the number of nodes increases, the density increases in both cost areas. 

A summary description of the instances is given in Table~\ref{table:instances}. The first column indicates the instance size ($|V|$). Then, the following columns indicate: the number of incompatible pairs (Pairs); the average number of potential access hubs per node (Avg Access); and, the average distance between any pair of nodes (Avg Dist). 

It can be observed that the number of incompatible pairs increases with the size of the instances. In these instances, as their size increases, the average number of accessible hubs per node also increases. The opposite happens with the average distance between pairs of nodes. In general, as the instance size increases, new nodes are added to the network in such a way that the average distance is reduced. This effect does not apply to instances with $|V|\in\{30,60,100\}$ whose average distance between pairs of nodes increases slightly. 

\begin{table}[ht]
\begin{center}
\caption{CAB instances detail.}
\label{table:instances}
\begin{tabular}{crrr}
\toprule
 & & \multicolumn{1}{l}{Avg} & \multicolumn{1}{l}{Avg} \\
$|V|$ & Pairs & \multicolumn{1}{l}{Access} & \multicolumn{1}{l}{Dist} \\ 
\midrule
10 & 9 & 2.20 & 1263 \\ 
20 & 14 & 2.80 & 1101 \\ 
30 & 27 & 4.67 & 1103 \\ 
40 & 39 & 7.35 & 1037 \\ 
50 & 55 & 7.52 & 1023 \\ 
\arrayrulecolor{black!30}\midrule\arrayrulecolor{black}
60 & 88 & 8.07 & 1045 \\ 
70 & 116 &  9.34 & 1042 \\ 
80 & 137 &  10.83 & 1016 \\ 
90 & 145 & 13.38 & 974 \\ 
100 & 182 & 14.98 & 973 \\ 
\bottomrule
\end{tabular}
\end{center}
\end{table}

\subsection{Comparison of the performance between formulations $F_1$ and $F_2$}

In order to assess the strength of the formulations we have compared the performance of $F_1$ and $F_2$ when applying the preprocessing procedures and the valid inequalities detailed in previous sections. For these experiments a time limit (TL) of $7,200$ s was set for each formulation and tested instance.

The effect of the preprocessing in $F_1$ and $F_2$ is shown in Table~\ref{table:F1-F2-sizes}. The first column indicates the instance size $|V|$. Then, for each formulation, the number of variables (Variables) and constraints (Constraints) are shown, followed by the percent reduction in the number of variables and constraints due to preprocessing, \% Red$_v$, and  \% Red$_c$, respectively.

It can be seen that the number of variables and constraints in $F_1$ is much larger than in formulation $F_2$. Even if the percent reductions achieved by preprocessing are higher for the first formulation, the resulting instance sizes are still much larger than those for the second. 

\begin{table}[ht]
\caption{Dimensions for $F_1$ and $F_2$ with and without preprocessing.}
\label{table:F1-F2-sizes}
\begin{center}
\begin{tabular}{crcrcrcrc}
\toprule
 & \multicolumn{4}{c}{$F_1$} 
 & \multicolumn{4}{c}{$F_2$}  \\ 
\cmidrule(r){2-5}\cmidrule(r){6-9}
$|V|$
& Variables & \% Red$_v$ & Constraints & \% Red$_c$ 
& Variables & \% Red$_v$ & Constraints & \% Red$_c$ \\ 
\midrule
10 & 9,110 & 66 & 5,681 & 55 & 1,110 & 30 & 1,911 & 52 \\
20 & 152,420 & 58 & 84,761 & 52 & 8,420 & 14 & 13,721 & 42 \\
30 & 783,930 & 56 & 420,241 & 52 & 27,930 & 11 & 44,431 & 41 \\
40 & 2,497,640 & 55 & 1,315,121 & 51 & 65,640 & 8 & 103,041 & 40 \\
50 & 6,127,550 & 54 & 3,192,401 & 51 & 127,550 & 7 & 198,551 & 39 \\
\arrayrulecolor{black!30}\midrule\arrayrulecolor{black}
60 & 12,747,660 & 54 & 6,595,081 & 51 & 219,660 & 7 & 339,961 & 39 \\
70 & 23,671,970 & 53 & 12,186,161 & 51 & 347,970 & 7 & 536,271 & 39 \\
80 & 40,454,480 & 53 & 20,748,641 & 51 & 518,480 & 6 & 796,481 & 39 \\
90 & 64,889,190 & 53 & 33,185,521 & 51 & 737,190 & 5 & 1,129,591 & 38 \\
100 & 99,010,100 & 52 & 50,519,801 & 51 & 1,010,100 & 5 & 1,544,601 & 38 \\
\bottomrule
\end{tabular}
\end{center}
\end{table}

The results obtained with the formulations with and without enhancements are summarized in Table~\ref{table:F1-F2-results}. The enhancements considered include the preprocessing procedures and the valid inequalities detailed in the previous sections. 
After some preliminary testing, we decided not to use inequalities~\eqref{vi_clique_Gen}, since they required a separation procedure that was not particularly effective.
In order to reduce the time required to begin the enumeration in the resolution of large instances, the maximum number of LPs solved by CPLEX at the root node was limited to five after some preliminary testing. 

The first column of the table, indicating the size of the instance, $|V|$, is followed by two blocks of columns, one for each formulation. In its turn, each block is divided into two sub-blocks, resulting in four variants in total. The first sub-block (Raw) shows the results obtained with CPLEX using the formulation as it is, while the second sub-block (Enhanced) shows the results obtained after applying the enhancements. Columns ``\% Gap'', ``T(s)'' and ``\% $\Delta$'', indicate the percent optimality gap, computing time (in seconds), and the percent deviation with respect to the optimal solution, respectively. The percent optimality gap is computed as $100\cdot(z - \underline z) / z$ where $z$ is the best integer objective value (upper bound) and $\underline z$ the best lower bound at termination. $\%\Delta=100\cdot(z - z^*) / z^*$ computes the percent deviation between the objective value of the best solution produced by the corresponding version ($z$) and the optimal value ($z^*$). For informational purposes, optimal solutions have been computed for all instances by running CPLEX with a larger time limit.

As the results show, formulation $F_1$ cannot obtain a feasible solution for any instance with $n\geq 50$ and only finds optimal solutions for small-size instances with $n\in\{10, 20\}$. The enhanced version of formulation $F_1$ finds a feasible solution to instances with up to $50$ nodes,  proving the optimality of the best solutions found for instances with up to $30$ nodes. 

On the other hand, formulation $F_2$ can obtain feasible solutions for all the instances considered, even without enhancements. The effects of the pre-processing procedures and valid inequalities in $F_2$ are evident in larger instances (more than 50 nodes), where the optimality gaps decrease and the quality of the feasible solutions improves.
Regarding the smaller instances (up to 50 nodes), solution time decreases in all cases except the 50-node instance, where the enhancements have the opposite effect. 

Formulation $F_2$ is able to prove optimality within the time limit for all the small instances (with and without enhancements) and, with enhancements, for the 70-node instance as well. It can be observed that the 90-node instance appears more difficult than the 100-node one, resulting in a larger optimality gap at termination.

These results confirm that formulation $F_2$ outperforms $F_1$; therefore, in the remainder of this section, we focus on the results obtained when applying the second formulation together with the proposed matheuristic. 

\begin{scriptsize}
\begin{table}[ht]
\caption{Comparison of $F_1$ and $F_2$ variants with a time limit of 7,200 s.}
\label{table:F1-F2-results}
\begin{center}
\begin{tabular}{crrrrrrrrrrrr}
\toprule
 & \multicolumn{6}{c}{$F_1$} 
 & \multicolumn{6}{c}{$F_2$}  \\ 
\cmidrule(r){2-7} \cmidrule(r){8-13} 
\multicolumn{1}{c}{ } 
 &  \multicolumn{3}{c}{Raw}  
 & \multicolumn{3}{c}{Enhanced} 
 & \multicolumn{3}{c}{Raw} 
 & \multicolumn{3}{c}{Enhanced} \\ 
\cmidrule(r){2-4} \cmidrule(r){5-7} \cmidrule(r){8-10} \cmidrule(r){11-13}
$|V|$
& \% Gap & T(s) & \%$\Delta$
& \% Gap & T(s) & \%$\Delta$
& \% Gap & T(s) & \%$\Delta$
& \% Gap & T(s) & \%$\Delta$\\ 
\midrule
10 & 0.00 & 1 & 0.00 & 0.00 & 0.11 & 0.00 & 0.00 & 0 & 0.00 & 0.00 & 0 & 0.00 \\
20 & 0.00 & 11 & 0.00 & 0.00 & 15.27 & 0.00 & 0.00 & 2 & 0.00 & 0.00 & 1 & 0.00\\
30 & 10.09 & TL & 6.08 & 0.00 & 650.23 & 0.00 & 0.00 & 23 & 0.00 & 0.00 & 9 & 0.00 \\
40 & 63.93 & TL & 137.04 & 3.63 & TL & 3.14 & 0.00 & 77 & 0.00 & 0.00 & 27 & 0.00 \\
50 & - & TL & - & 96.15 & TL & 287.25 & 0.00 & 768 & 0.00 & 0.00 & 1972 & 0.00\\
\arrayrulecolor{black!30}\midrule\arrayrulecolor{black}
60 & - & TL & - & - & TL & - & 5.57 & TL & 2.85 & 2.92 & TL & 0.87 \\
70 & - & TL & - & - & TL & - & 30.93 & TL & 39.11 & 0.00 & 1988 & 0.00 \\
80 & - & TL & - & - & TL & - & 55.38 & TL & 84.48 & 5.41 & TL & 0.97 \\
90 & - & TL & - & - & TL & - & 72.59 & TL & 186.79 & 26.18 & TL & 16.52 \\
100 & - & TL & - & - & TL & - & 75.80 & TL & 223.04 & 13.41 & TL & 8.04 \\
\bottomrule
\end{tabular}
\end{center}
\scriptsize
{\bf Raw}: original formulation.\\
{\bf Enhanced}: formulation with enhancements (pre-processing, valid inequality and limited LP-resolutions).\\
\end{table}
\end{scriptsize}

\subsection{Matheuristic Performance}

In this section, we analyze the performance of the matheuristic. Before starting, it is worth mentioning that almost all the solutions obtained are optimal. The matheuristic produces near-optimal solutions only for the instances with 80 and 90 nodes.

Table~\ref{table:MH} summarizes the performance of the matheuristic with parameters $\rho_{max}=5$ (maximum neighborhood radius), $\eta_{max}=|V|/5$ (maximum number of explored solutions in the same neighborhood), and $\tau_{max}=60$ s. (maximum computing time).

The first column indicates the instance size, $|V|$. 
Then, for each improved solution produced by the procedure, the iteration and neighborhood radius ($t$, $\rho$) at which it was obtained are reported.
Since, in all instances, the solution from the first iteration is obviously the best at that stage of the procedure, it is not included in the table. 
The last four columns stand for the iteration producing the best solution ($t^*$), the total number of iterations ($t$), the average computing time per iteration ($\overline\tau$), and the total computing time ($\tau$). 

In almost all instances, the average time per iteration increases with the size of the instance, reflecting that the difficulty of the subproblems increases. Only in the smallest
instance, the procedure terminates because the location subproblem is infeasible, after exploring all the neighborhoods. 

Instances with 90 and 100 nodes reach the time limit; therefore, the algorithm does not have enough time to fully explore the solutions' neighborhoods. 

\begin{table}[hbt]
\begin{center}
\caption{Matheuristic performance for $\rho_{max}=5$, $\eta_{max}=|V|/5$ and $\tau_{max}=60$ s.}
\label{table:MH}
\begin{tabular}{clllllcccr}
\toprule
$|V|$ 
& \multicolumn{5}{c}{\text{Steps with improved solutions} $(t,\rho)$}
& $t^*$ & $t$ & $\overline \tau$ & \multicolumn{1}{c}{$\tau$} \\ 
\midrule
10 & (2,1)    & (3,1)  & (6,2)  &  		& 		& 6 & 11 & 0.00 & 0.04 \\ 
20 & (2,1)    &        &        &  		& 		& 2 & 20 & 0.01 & 0.25 \\ 
30 & (2,1)    & (3,1)  & (21,3) &  		& 		& 21 & 49 & 0.03 & 1.64 \\ 
40 & (2,1)    &        & (15,2) & (16,1)& 		& 16 & 54 & 0.03 & 1.59 \\ 
50 & (2,1)    & (3,1)  & (46,5) &  		& 		& 46 & 94 & 0.28 & 26.22 \\ 
\arrayrulecolor{black!30}\midrule\arrayrulecolor{black}
60 & (2,1)    &        & (19,2) & (62,4)& (84,2)& 84 & 142 & 0.36 & 50.98 \\ 
70 & (2,1)    &        & (19,2) &  		& 		& 19 & 87 & 0.31 & 27.40 \\ 
80 & (2,1)    &        &        &  		& 		& 2 & 80 & 0.22 & 17.26 \\ 
90 & (2,1)    &        & (21,2) &  		& 		& 21 & 102 & 0.59 & TL \\ 
100&          &        & (23,2) & (44,2) & 		& 44 & 134 & 0.45 & TL \\ 
\bottomrule
\end{tabular}
\end{center}
\end{table}

\subsection{Performance of $F_2$ with enhancements}

We next analyze the results of three different variants of formulation $F_2$, which is the one performing best.
The four variants considered are the following: original formulation (``Raw''), formulation with enhancements (``Enhanced''), and, formulation with enhancements using as starting solution the matheuristic solution (``Warm Start''). In this section, a time limit (TL) of 14,400 s was set for each tested instance.

Table~\ref{table:F2:cplex} summarizes the results produced by CPLEX for all three variants.
For each of them, the headers ``\#N'', ``\%Gap'', ``T(s)'', and ``\%$\Delta$'', indicate the number of explored nodes, the optimality gap (in percentage), the computing time (in seconds), and percent optimality deviation of the best solution at termination, respectively. 
To avoid repetition, results are shown only for large instances (more than 50 nodes), since all smaller instances (up to 50 nodes) were already solved to optimality within a time limit of 7,200 s.

As shown in Table~\ref{table:F1-F2-results}, across all instances, incorporating the warm start reduces computation time for both the raw and enhanced versions. In particular, the 50-node instance is solved in $379$ seconds when using the matheuristic solution as a warm start. 

In terms of optimality gaps and computation times, the enhanced version with warm start achieves the largest reductions in both cases, with optimality gaps below 6\% for all large instances. Note that the solution provided by the matheuristic and used as a warm start is not improved during the branching, so the decrease in the optimality gap is only due to an improvement of the upper bound. 

As can be seen in the table, the solution provided by the matheuristic is optimal in all cases except for the 80 and 90 node instances, where the solution attained has a deviation with respect to the optimal solution of $0.33\%$ and $0.31\%$ respectively.

\begin{table}[htb]
\caption{Performance of $F_2$ variants with a time limit of $14\,400$ s.}
\label{table:F2:cplex}
\begin{center}
\begin{tabular}{crrrrrrrrrrrr}
\toprule
& 
\multicolumn{4}{c}{Raw} & 
\multicolumn{4}{c}{Enhanced} & 
\multicolumn{4}{c}{Warm Start} \\ 
\cmidrule(r){2-5}\cmidrule(r){6-9}\cmidrule(r){10-13}
$|V|$
& \multicolumn{1}{c}{\#N} & \multicolumn{1}{c}{\%Gap} &  \multicolumn{1}{c}{T(s)} & \multicolumn{1}{c}{\%$\Delta$}
& \multicolumn{1}{c}{\#N} & \multicolumn{1}{c}{\%Gap} &  \multicolumn{1}{c}{T(s)} & \multicolumn{1}{c}{\%$\Delta$}
& \multicolumn{1}{c}{\#N} & \multicolumn{1}{c}{\%Gap} &  \multicolumn{1}{c}{T(s)} & \multicolumn{1}{c}{\%$\Delta$}\\ 
\midrule
60 & 1345 & 2.70 & TL & 0.91 & 2490 & 2.50 & TL & 0.87 & 2583 & 1.71 & TL & 0.00 \\ 
70 & 833 & 16.81 & TL & 15.51 & 1181 & 0.00 & 1933 & 0.00 & 353 & 0.00 & 1161 & 0.00 \\ 
80 & 1110 & 8.58 & TL & 3.77 & 1065 & 5.41 & TL & 0.97 & 1208 & 3.96 & TL & 0.33 \\ 
90 & 1 & 71.05 & TL & 172.44 & 466 & 25.08 & TL & 16.52 & 1001 & 5.77 & TL & 0.31 \\ 
100 & 155 & 33.29 & TL & 22.03 & 459 & 12.47 & TL & 6.87 & 1377 & 3.73 & TL & 0.00 \\ 
\bottomrule
\end{tabular}
\end{center}
\scriptsize
{\bf Raw}: original formulation.\\
{\bf Enhanced}: formulation with enhancements (pre-processing, valid inequality and limited LP-iterations at the root node).\\
{\bf Warm Start}: enhanced formulation with warm start. 
\end{table}

\subsection{Managerial insights}

In this section, some insights into the structure of the best-known solutions are provided. 
Table~\ref{table:best-sols} summarizes some details of the obtained solutions. The first two columns indicate the instance size, $|V|$, and the objective value attained, $z$, respectively. 
Then, the next five columns provide the following information about the resulting network: 
number of hubs (H), number of uncovered nodes (U), percentage of covered nodes (\%C), interhub network density (Dens), and average interhub travel time (Avg Time). 
Finally, the last two columns refer to covered nodes only, indicating the 
average number of activated access hubs (Avg Act Access),
and average distance to the activated access hub (Avg Dist).
The interhub network density is computed as $2\, M / (H \, (H-1))$ where $H$ denotes the number of hubs and $M$ the number of interhub edges in the network. This value is one when the interhub network is a complete graph. 

\begin{table}[htb]
\begin{center}
\caption{Some characteristics of best-known solutions.}
\label{table:best-sols}
\begin{tabular}{crrrrrccc}
\toprule
 & & & & & & \multicolumn{1}{l}{Avg} & \multicolumn{1}{l}{Avg Act } & \multicolumn{1}{l}{Avg}\\
$|V|$ & 
\multicolumn{1}{c}{$z$} & 
H & 
U & 
\%C &
Dens & 
\multicolumn{1}{l}{Time} & 
\multicolumn{1}{l}{Access} & 
\multicolumn{1}{l}{Dist }
\\ \midrule
10 & 1,899,515 & 2 & 4 & 60 & 1.00 & 48.00 & 1.00 & 121.33 \\ 
20 & 3,587,403 & 4 & 7 & 65 & 0.50 & 86.50 & 1.15 & 163.31 \\ 
30 & 3,726,744 & 6 & 5 & 83 & 0.33 & 150.47 & 1.20 & 158.04 \\ 
40 & 3,817,961 & 6 & 5 & 88 & 0.33 & 134.33 & 1.26 & 186.46 \\ 
50 & 4,059,413 & 8 & 4 & 92 & 0.25 & 130.54 & 1.04 & 187.07 \\ 
\arrayrulecolor{black!30}\midrule\arrayrulecolor{black}
60 & 4,574,324 & 9 & 5 & 92 & 0.25 & 134.08 & 1.15 & 177.33 \\ 
70 & 4,234,964 & 9 & 3 & 96 & 0.22 & 168.89 & 1.13 & 196.45 \\ 
80 & 4,393,883 & 9 & 2 & 98 & 0.31 & 151.36 & 1.15 & 201.24 \\ 
90 & 4,666,686 & 10 & 2 & 98 & 0.27 & 153.64 & 1.27 & 195.97 \\ 
100 & 4,502,102 & 10 & 2 & 98 & 0.22 & 156.38 & 1.30 & 189.89 \\ 
\bottomrule
\end{tabular}
\end{center}
\end{table}

\begin{figure}[htb]
\begin{subfigure}[b]{0.35\textwidth}
\frame{\includegraphics[width=\textwidth]{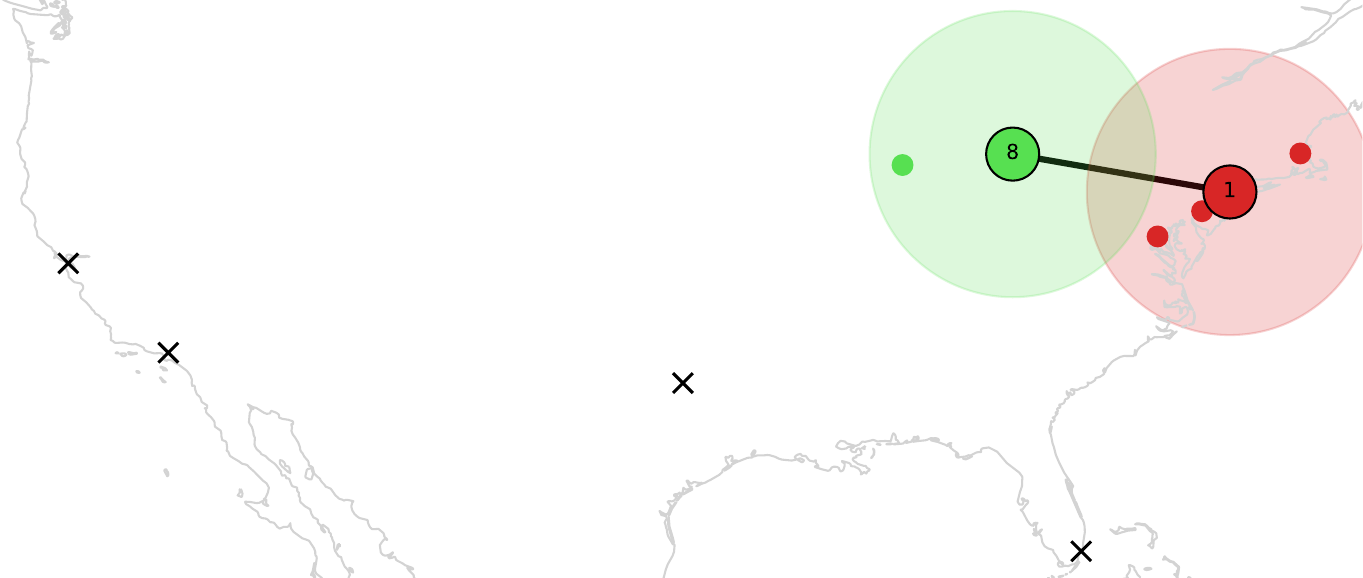}}
\caption{$|V|=10$}
\end{subfigure}
\begin{subfigure}[b]{0.35\textwidth}
\frame{\includegraphics[width=\textwidth]{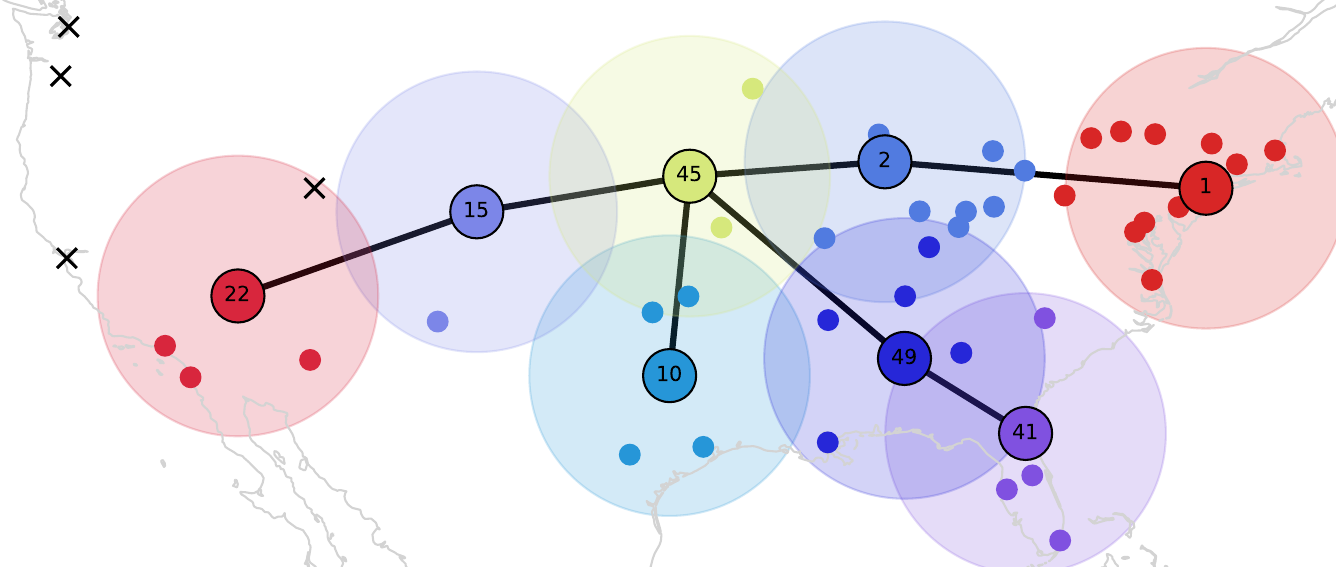}}
\caption{$|V|=50$}
\end{subfigure}
\begin{subfigure}[b]{0.35\textwidth}
\frame{\includegraphics[width=\textwidth]{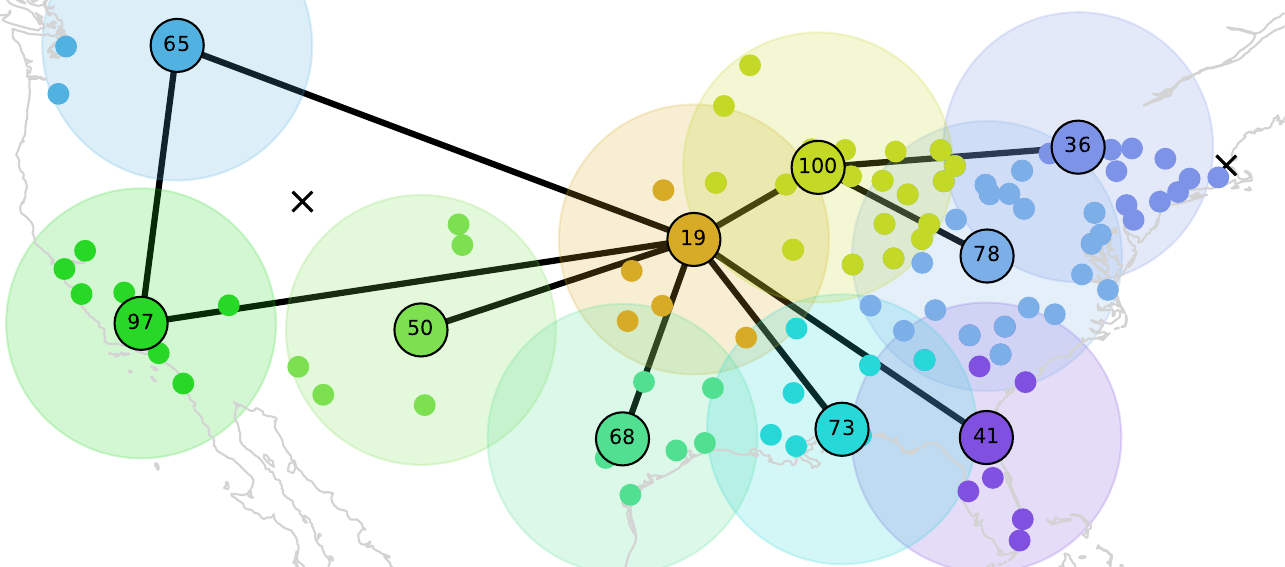}}
\caption{$|V|=100$}
\end{subfigure}
\caption{Maps of optimal solutions for $T_{max}=230$ and $\delta=340$.}
\label{fig:best-sols:mapas}
\end{figure}

The density of the interhub network remains below $0.4$ for the largest instances. As the size of the instances increases, and solutions can cover almost all nodes with a lower percentage of activated hubs. The objective value attains its peak around $4.6$ million monetary units (for 60 and 90 node instances), and all other cases remain below this value. Average travel times between two hubs are around $146$ for instances with 30 nodes onward. In all instances, the average distance of covered nodes to their access hubs is around $177$, with a hub coverage radius of $\delta=340$. The average number of activated access hubs per covered node does not exceed $1.3$; that is, covered nodes are, on average, in the coverage area of one or two hubs. This is illustrated in Figure~\ref{fig:best-sols:mapas}, which shows some of the solutions.

\begin{figure}[htb]
\hspace*{0.50cm}
\includegraphics[height=4cm]{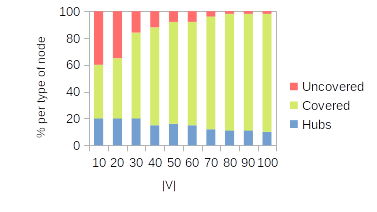}\hspace{0.5cm}
\hspace*{0.50cm}
\includegraphics[height=4cm]{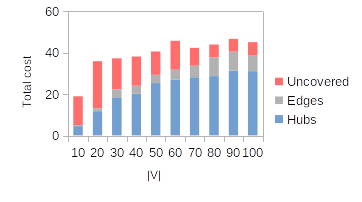}
\vspace{-0.50cm}
\caption{Left: percent distribution per type of node. Right: contribution to the total cost of each objective function term. } 
\label{fig:best-sols:nodes-cost}
\end{figure}

Figure~\ref{fig:best-sols:nodes-cost}~(left) shows, for each instance, the percentage of open hubs, covered nodes, and uncovered nodes in its best-known solution. Due to the structure of the instances, the percentage of uncovered nodes decreases with increasing instance size. Note that, for large instances, this percentage remains below $10\%$. For large instances, the number of open hubs is nearly constant (9 or 10 hubs), and therefore, the percentage of installed hubs decreases as instance size increases. This percentage stays between 10\% and 20 \% for all the instances.

Figure~\ref{fig:best-sols:nodes-cost}~(right) illustrates the contribution to the total cost of ($i$) the penalty due to uncovered nodes (in red), ($ii$) the setup costs of activated hubs (in blue), and ($iii$) the setup costs of activated interhub edges (in gray). Hub installation costs remain nearly constant for large instances because the number of hubs is nearly the same, even though their locations change slightly, resulting in small variations in activation costs. As mentioned, the percentage of uncovered nodes decreases with increasing instance size. Therefore, a similar behavior can be observed with respect to the penalty costs. The setup costs of interhub edges do not increase linearly with the instance size. Note that this cost not only depends on the number of edges installed but also on the length of each activated edge. Optimal solutions for small instances have between 1 and 7 edges, whereas for large instances it ranges from 8 to 12. 

Figure~\ref{fig:best-sols:mapas} depicts the maps corresponding to optimal solutions for the instances with 10, 50 and 100 nodes. Solution maps for all the instances are given in the Appendix. 

Note that the optimal networks for the 10- and 50-node instances correspond to minimum-cost spanning trees in the graphs induced by the activated hubs. That is, the activated interhub edges are those that minimize the edge setup costs, provided that travel times do not exceed $T_{max}$. However, in the 100-node optimal network, an extra edge (connecting the hubs located at nodes 65 and 97) is added to attain feasibility, creating a cycle. 

\subsection{Sensitivity analysis}

In this section, we analyze the effect on the optimal solutions of the two parameters considered: coverage radius, $\delta$, and maximum transit network travel time, $T_{max}$. We perform a sensitivity analysis using the 50-node instance as a test case. The parameters default values used were $T_{max}=230$ and $\delta=340$, and we  explore solutions with $\delta$ varying  in $\{240,290,340,390,440\}$ and  $T_{max}$ varying in $\{130, 180, 230, 280, 330\}$.

Table~\ref{table:sensitivity:solutions:Tmax} describes features of the instances and optimal solutions for the different values of $T_{max}$, whereas 
Table~\ref{table:sensitivity:solutions:delta} gives the same information when varying $\delta$.
The characteristics of the instance are detailed in the first four columns of the tables: 
maximum travel time ($T_{max}$), number of incompatible pairs (Pairs), and average number of access hubs per node (Avg Access).
The remaining columns describe the network details as follows: optimal objective value ($z^*$), number of hubs (H), number of uncovered nodes (U), percentage of covered nodes (\%C), interhub network density (Dens), average travel time between two hubs (Avg Time), average number of activated access hubs per covered node (Avg Act Access), and average distance to an access hub per covered node (Avg Dist).
For reference, the results corresponding to the default configuration are highlighted in gray. 

\begin{figure}[htb]
\begin{center}
\includegraphics[height=4.5cm]{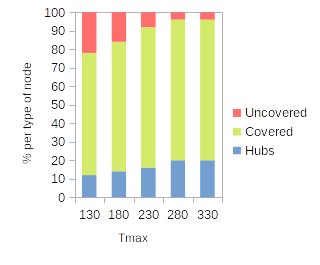}
\hspace*{1cm}
\includegraphics[height=4.5cm]{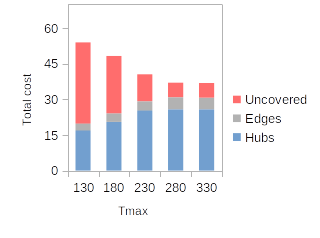}
\vspace{-0.50cm}
\caption{(Left) Percent distribution per type of node; (Right) Contribution to the objective function of each of its three terms; Results with $\delta=340$ and $T_{max}\in \{130, 180, 230, 280, 330\}$.}
\label{fig:sensitivity:Tmax}
\end{center}
\end{figure}

\begin{figure}[htb]
\begin{subfigure}[b]{0.35\textwidth}
\frame{\includegraphics[width=\textwidth]{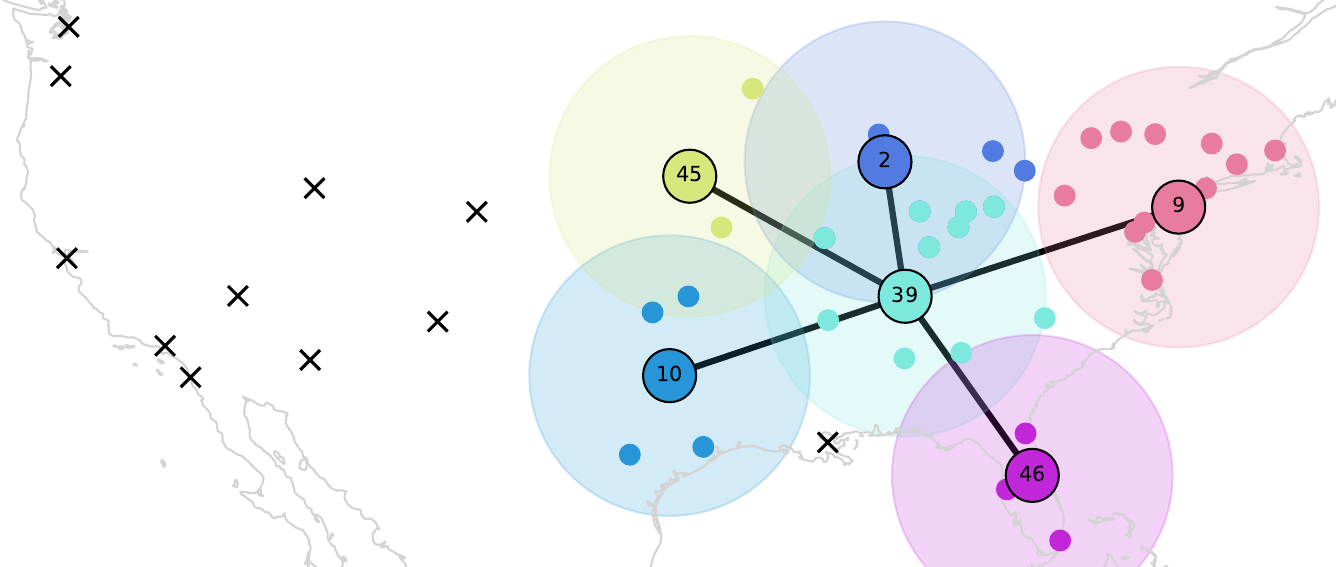}}
\caption{$\delta=340$ and $T_{max}=130$}
\end{subfigure}
\begin{subfigure}[b]{0.35\textwidth}
\frame{\includegraphics[width=\textwidth]{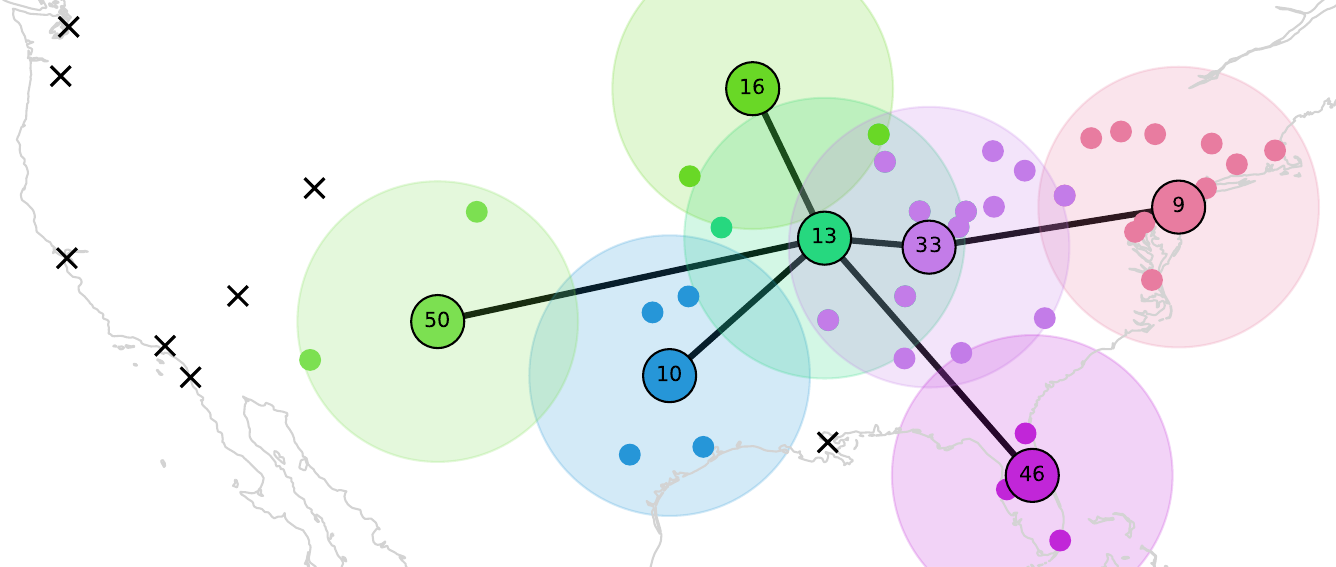}}
\caption{$\delta=340$ and $T_{max}=180$}
\end{subfigure}
\begin{subfigure}[b]{0.35\textwidth}
\frame{\includegraphics[width=\textwidth]{50CAB_R340_T230_network.pdf}}
\caption{$\delta=340$ and $T_{max}=230$ (default)}
\end{subfigure}
\\
\begin{subfigure}[b]{0.35\textwidth}
\frame{\includegraphics[width=\textwidth]{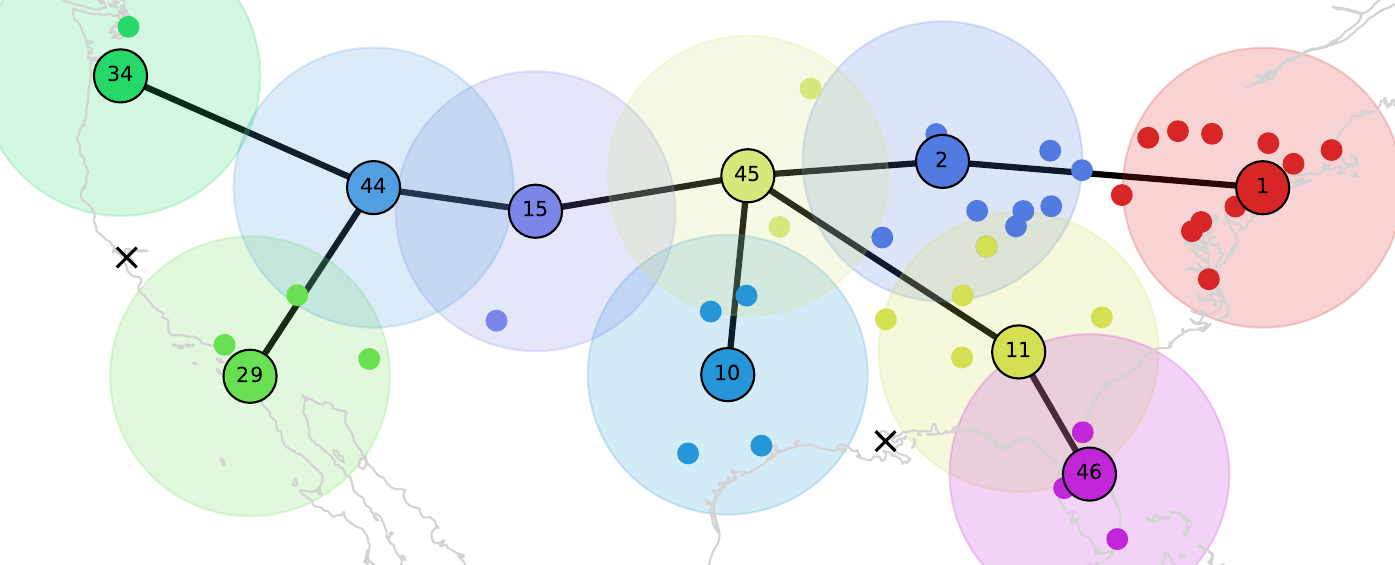}}
\caption{$\delta=340$ and $T_{max}=280$}
\end{subfigure}
\begin{subfigure}[b]{0.35\textwidth}
\frame{\includegraphics[width=\textwidth]{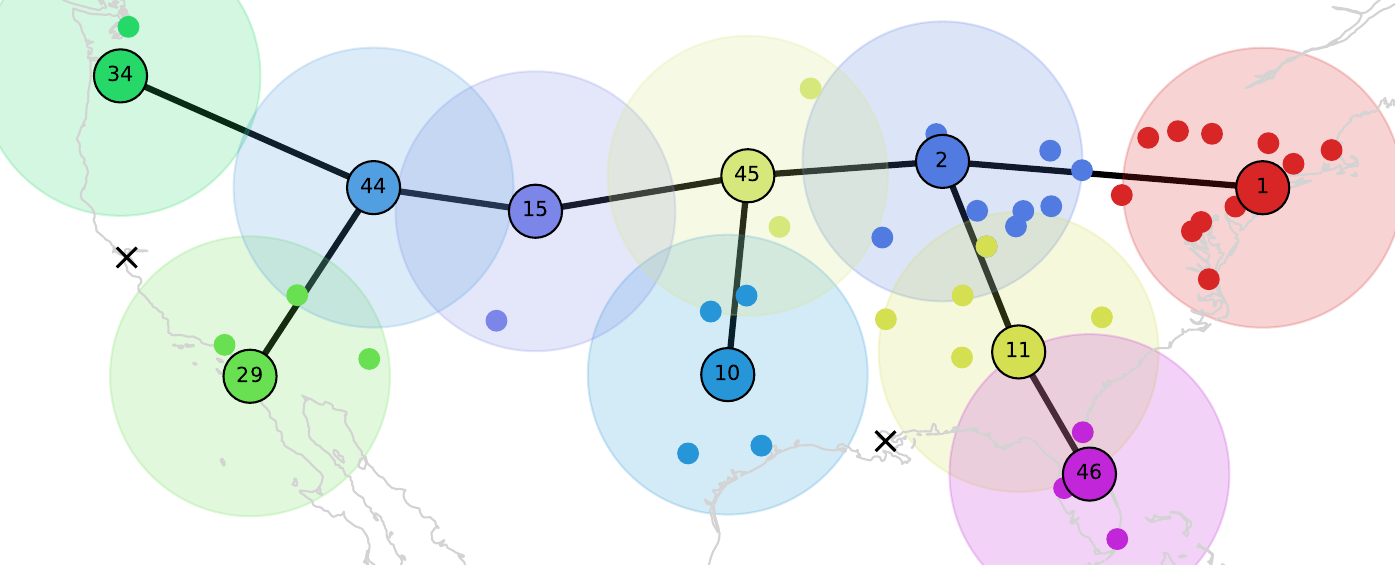}}
\caption{$\delta=340$ and $T_{max}=330$}
\end{subfigure}
\caption{Optimal solutions for instance with $|V|=50$, $\delta=340$ and $T_{max}\in \{130, 180, 230, 280, 330\}$.}
\label{fig:sensitivity:mapas:Tmax}
\end{figure}

\begin{figure}[htb]
\begin{center}
\includegraphics[height=4.3cm]{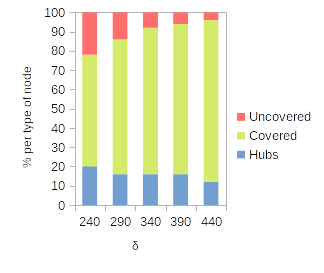}
\hspace*{1cm}
\includegraphics[height=4.3cm]{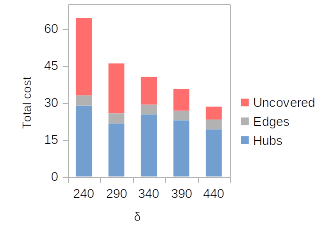}
\vspace{-0.50cm}
\caption{(Left) Percent distribution per type of node; (Right) Contribution to the objective function of each of its three terms; Results with $\delta=230$, $T_{max}\in \{240,290,340,390,440\}$.}
\label{fig:sensitivity:delta}
\end{center}
\end{figure}

\begin{figure}[htb]
\begin{subfigure}[b]{0.35\textwidth}
\frame{\includegraphics[width=\textwidth]{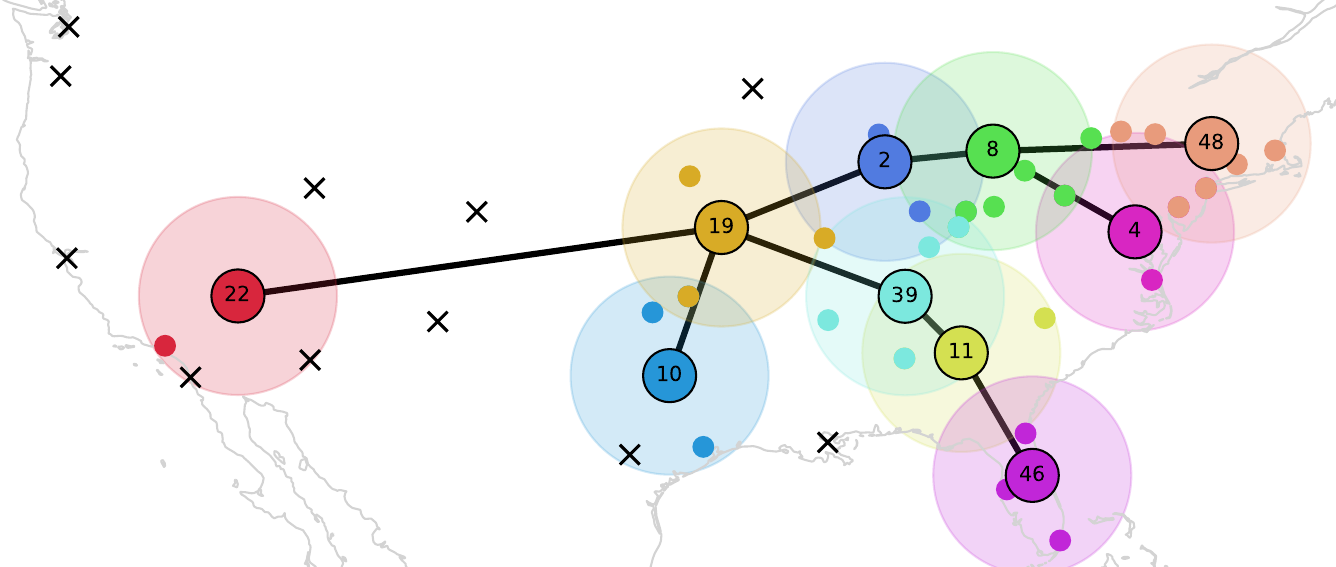}}
\caption{$\delta=240$ and $T_{max}=230$}
\end{subfigure}
\begin{subfigure}[b]{0.35\textwidth}
\frame{\includegraphics[width=\textwidth]{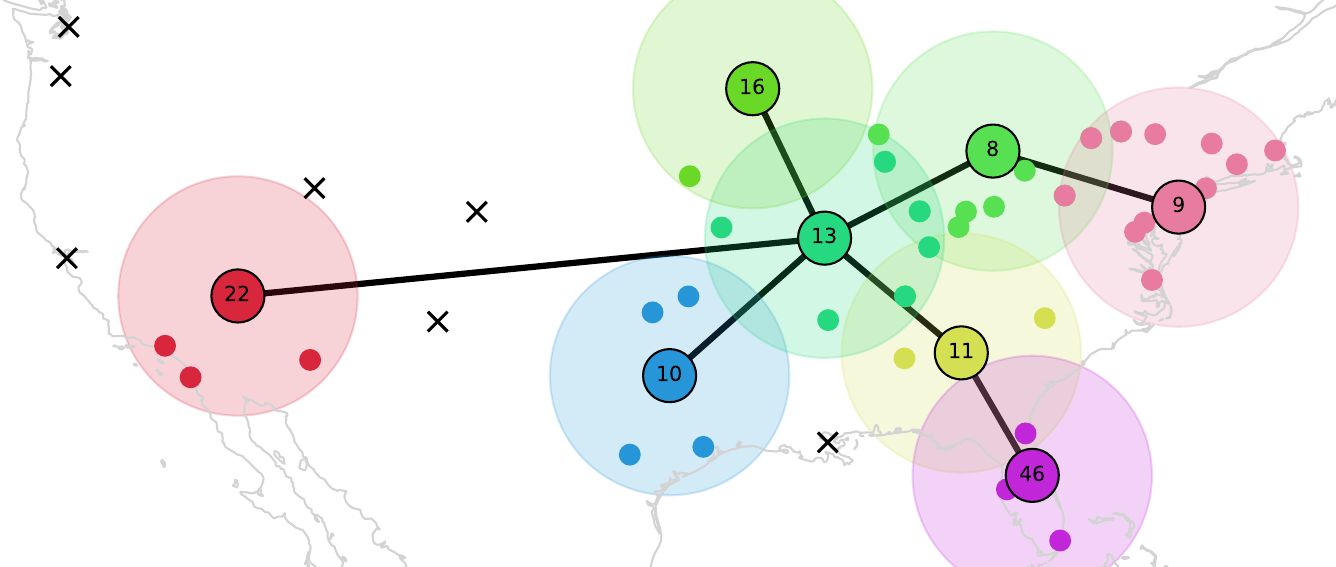}}
\caption{$\delta=290$ and $T_{max}=230$}
\end{subfigure}
\begin{subfigure}[b]{0.35\textwidth}
\frame{\includegraphics[width=\textwidth]{50CAB_R340_T230_network.pdf}}
\caption{$\delta=340$ and $T_{max}=230$ (default)}
\end{subfigure}
\\
\begin{subfigure}[b]{0.35\textwidth}
\frame{\includegraphics[width=\textwidth]{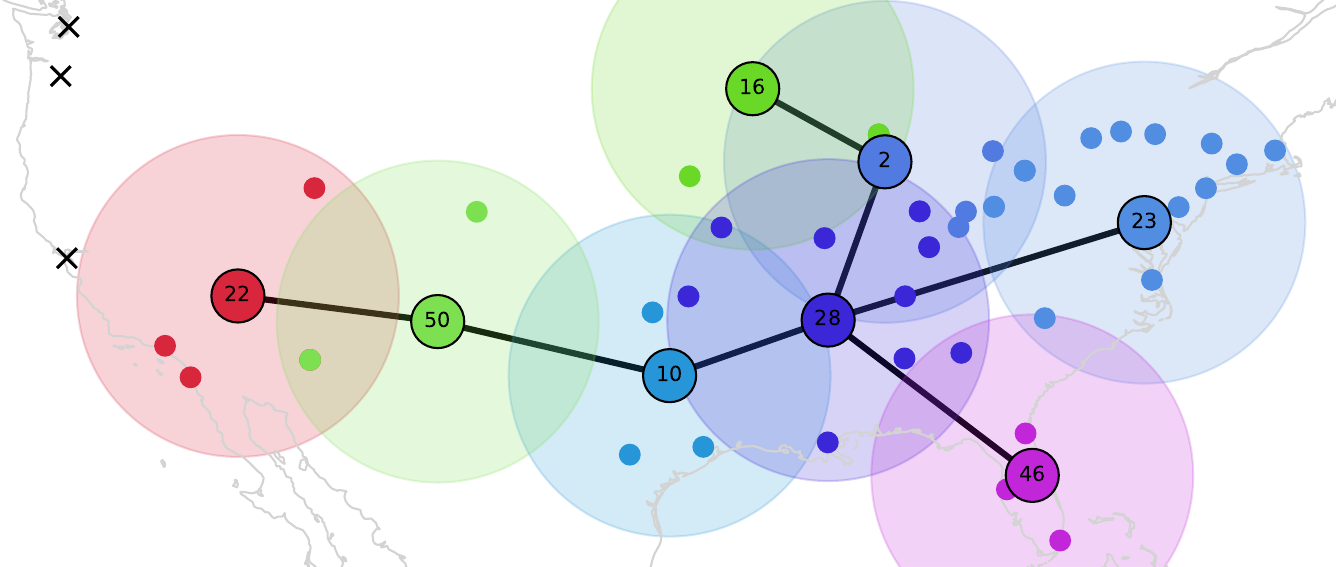}}
\caption{$\delta=390$ and $T_{max}=230$}
\end{subfigure}
\begin{subfigure}[b]{0.35\textwidth}
\frame{\includegraphics[width=\textwidth]{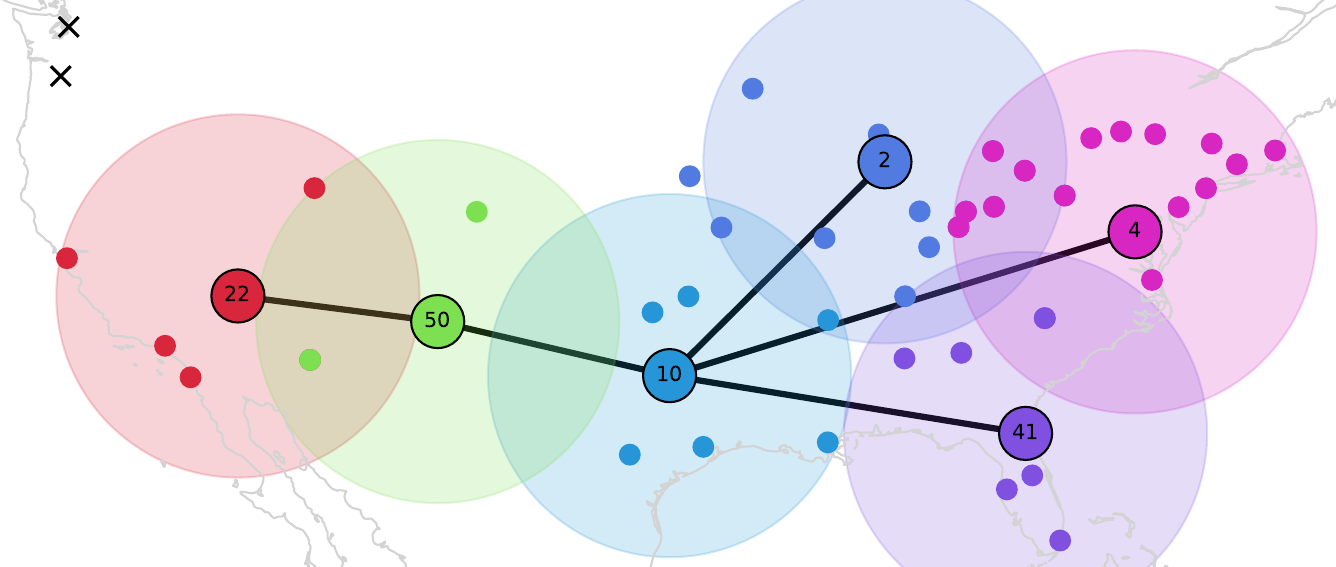}}
\caption{$\delta=440$ and $T_{max}=230$}
\end{subfigure}
\caption{Optimal solutions for instance $|V|=50$ with $T_{max}=230$ and $\delta\in \{240,290,340,390,440\}$.}
\label{fig:sensitivity:mapas:delta}
\end{figure}

Figure~\ref{fig:sensitivity:Tmax} complements the above information by showing 
the percent distribution per type of node (left), as well as the contribution to the objective function value of each of its three terms (right) when varying $T_{max}$. Then, the optimal solutions for each $T_{max}$ value are depicted in Figure~\ref{fig:sensitivity:mapas:Tmax}. The same information is illustrated in Figures~\ref{fig:sensitivity:delta} and~\ref{fig:sensitivity:mapas:delta} under variations of the $\delta$ parameter.

As expected, relaxing the constraints on the maximum travel time within the transit network or the coverage radius, by increasing the values of $T_{max}$ or $\delta$, respectively, yields solutions with lower objective function values, $z^*$. As shown in the bar charts, the reduction in total cost is due to a decrease in the penalty term for uncovered nodes. 

\begin{table}[htb]
\begin{center}
\caption{Instances and networks details with $\delta=340$ and $T_{max}\in \{130, 180, 230, 280, 330\}$.}
\label{table:sensitivity:solutions:Tmax}
\begin{tabular}{cccccccccccc}
\toprule
\multicolumn{4}{c}{Instance} & \multicolumn{8}{c}{Optimal Solution}\\
\cmidrule(lr){1-4}\cmidrule(lr){5-12}
 & & &  \multicolumn{1}{l}{Avg} & & & & & &  \multicolumn{1}{l}{Avg} & Avg Act &  \multicolumn{1}{l}{Avg} \\
$\delta$ & $T_{max}$ & Pairs & Access 
& $z^*$ & H & U & \%C  & Dens & \multicolumn{1}{l}{Time} &  \multicolumn{1}{l}{Access} & \multicolumn{1}{l}{Dist}
\\ \midrule
 & 130 & 348 & & 5,402,008 & 6 & 11 & 78 & 0.33 & 101.67 & 1.15 & 175.79 \\ 
 & 180 & 186 & & 4,829,879 & 7 & 8 & 84 & 0.29 & 122.62 & 1.24 & 186.26 \\ 
 \rowcolor{shadecolor} 
340 & 230 & 55 & 7.52 & 4,059,413 & 8 & 4 & 92 & 0.25 & 130.54 & 1.04 & 187.07 \\ 
 & 280 & 0 & & 3,710,941 & 10 & 2 & 96 & 0.20 & 157.56 & 1.04 & 169.94 \\ 
 & 330 & 0 & & 3,688,541 & 10 & 2 & 96 & 0.20 & 159.07 & 1.04 & 169.94 \\ 
\bottomrule
\end{tabular}
\end{center}
\end{table}

It can also be observed in the table that as $T_{max}$ increases, both the number of installed hubs and the percentage of covered nodes increase, while the number of uncovered nodes decreases. As can be seen in the figures, the configuration of the interhub network expands from a star spanning hubs near the east coast to a tree with a larger diameter, spanning hubs from both coastal areas. Indeed, as the diameter of the transit network increases, the average access distance per covered node decreases.

Note also that, as shown in the table, the number of incompatible pairs decreases as $T_{max}$ increases. This decrease also reduces the difficulty of the problem. Moreover, when there are no incompatible pairs, the optimal transit network corresponds to the minimum cost spanning tree in the graph induced by the activated hubs, relative to the edges activation costs.

In summary, allowing longer travel times in the transit network produces more \textit{expanded} hub networks with lower density. On the contrary, increasing the coverage radius has the opposite effect, as we discuss below.

\begin{table}[htb]
\begin{center}
\begin{small}
\caption{Instances and networks details with $T_{max}=230$ and $\delta\in \{240,290,340,390,440\}$.}
\label{table:sensitivity:solutions:delta}
\begin{tabular}{cccc|cccccccc}
\toprule
\multicolumn{4}{c}{Instance} & \multicolumn{8}{c}{Optimal Solution}\\
\cmidrule(lr){1-4}\cmidrule(lr){5-12}
 & & &  \multicolumn{1}{l|}{Avg} & & & & & &  \multicolumn{1}{l}{Avg} & Avg Act &  \multicolumn{1}{l}{Avg} \\
$\delta$ & $T_{max}$ & Pairs & Access 
& $z^*$ & H & U & \%C  & Dens & \multicolumn{1}{l}{Time} &  \multicolumn{1}{l}{Access} & \multicolumn{1}{l}{Dist}
\\ \midrule
240 &  &  &  4.32 & 6,444,027 & 10 & 11 & 78 & 0.20 & 129.51 & 1.28 & 120.15 \\ 
290 &  &  &  5.76 & 4,608,569 & 8 & 7 & 86 & 0.25 & 137.21 & 1.16 & 165.79 \\ 
 \rowcolor{shadecolor} 
340 & 230 & 55 & 7.52 & 4,059,413 & 8 & 4 & 92 & 0.25 & 130.54 & 1.04 & 187.07 \\ 
390 &  &  & 8.64 &  3,562,091 & 8 & 3 & 94 & 0.25 & 130.93 & 1.21 & 200.19 \\ 
440 &  &  & 10.56 & 2,851,119 & 6 & 2 & 96 & 0.33 & 160.87 & 1.17 & 239.83 \\ 
\bottomrule
\end{tabular}
\end{small}
\end{center}
\end{table}

\vspace{-0.50cm}

Focusing now on the effect of the coverage radius parameter in the instance configuration, the average number of potential access hubs per node increases with the coverage radius $\delta$. In turn, this decreases the number of activated hubs in the optimal solutions and the penalty term in the objective function due to uncovered nodes, thereby reducing the optimal objective function value. 
As shown in Figure~\ref{fig:sensitivity:mapas:delta}, all solutions obtained with $T_{max}=230$ and varying coverage radii $\delta$ activate hubs in both the east and west coastal areas.
As $\delta$ increases, the table shows that the number of activated hubs decreases, whereas the percentage of covered nodes increases. 
As expected, because activated hubs remain located in both the west and east areas, the average travel time between a pair of activated hubs increases up to approximately 160. 

\subsection{Public transport examples}

In this section, we apply formulation F2 to a set of instances adapted from~\cite{Delle2025}, which emulate the use of public transport networks, specifically bus lines, for last-mile delivery (LMD). We use these public transport networks as realistic cases to illustrate our model. 

The LMD instances of~\cite{Delle2025} describe a network with two disjoint sets: the set of demand nodes ($D$) and the set of public transport stops ($S$). The union of these two sets forms the set of nodes of the input network $N=(V, A)$, where $V=D\cup S$ and $A=\{(i, j): i, j\in V\}$. We consider two instances with a total of $50$ demand nodes and the number of stops varying in $\{28, 31\}$. After some preliminary testing, we set parameters $\delta=20$ and $T_{max}=100$. 
For each node $i\in V$, the original LMD instances include the following data in addition to the type of node: activation cost ($g_i$), and coordinates (latitude and longitude). From the latitude and longitude, we compute the Euclidean distance between each pair of nodes denoted as $d_{ij}$ and the travel times $t_{ij} = \tau \, d_{ij}$, where $\tau = 0.2$ is the bus speed distance/time conversion factor. 

In the LMD instances, the setup costs $g_i$ are zero for demand nodes and non-zero for stop nodes. 
We consider that hubs can be located only at stop nodes, excluding demand nodes from the transportation network. That is, we set $G_i= 1\,000\, g_i$ for all $i\in S$ and $G_i=+\infty$ for all $i\in D$. 
Therefore, in the formulation, variables $y_i$ will be zero for all $i\in D$. 
The activation costs for interhub edges are set as $F_{km} = \sigma\,d_{km}$ for each $km\in E_H$ with $\sigma = 100$ a conversion factor.
Finally, we assume that all demand nodes must be covered by the resulting network, whereas stop nodes need not have access to the network (as they do not have demand). That is, we set the penalty values as follows: $P_i=G_i$ for all $i\in S$ and $P_i=+\infty$ for all $i\in D$.
In the formulation, variables $v_i$ will be zero for all demand nodes.

\begin{figure}[htb]
\begin{center}
\begin{subfigure}[b]{0.30\textwidth}
\frame{\includegraphics[width=\textwidth]{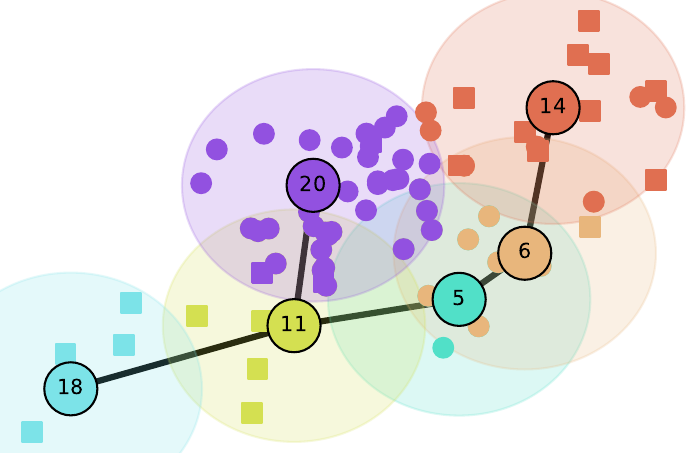}}
\caption{Instance with $|D|=50$ and $|S|=28$}
\end{subfigure}
\hspace{0.5cm}
\begin{subfigure}[b]{0.30\textwidth}
\frame{\includegraphics[width=\textwidth]{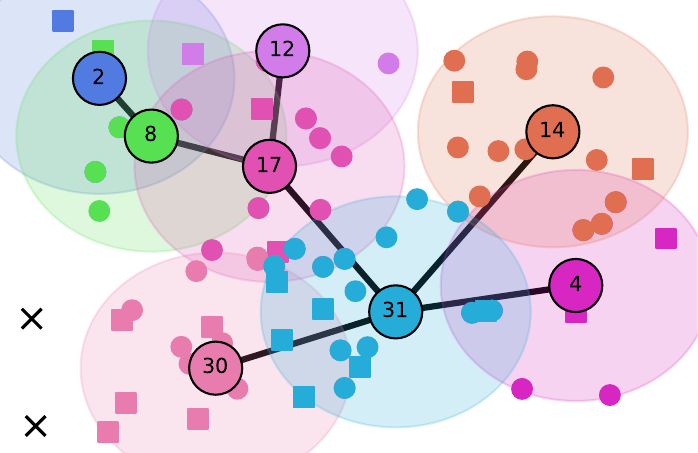}}
\caption{Instance with $|D|=50$ and $|S|=31$}
\end{subfigure}
\caption{Optimal solutions for two public transport instances.}
\label{fig:bus}
\end{center}
\end{figure}
\vspace{-0.50cm}

The optimal solutions for these instances are shown in Figure~\ref{fig:bus}, where stop nodes are represented by squares and demand nodes by circles. Bus stops that are neither activated nor covered are represented by a cross. 
As can be seen, the demand nodes in the instance are concentrated in the center of the urban area, while the bus line stops appear to extend towards the periphery, where demand nodes are less dense. 

In the solutions for the two instances, the percentages of activated bus stops are 21\% and 26\%, respectively. We can see that many stops on the periphery are not activated as part of the network, even though they are covered by it. In Figure \ref{fig:bus}(a), all nodes are covered by the resulting network, whereas in Figure \ref{fig:bus}(b), two stops are neither activated as part of the hub network nor covered by it. These two stops are located in the bottom-left part of the figure, somewhat farther from the other nodes and without demand nodes nearby.

\section{Conclusions}
\label{sec:conclu}

In this paper, we propose a modeling approach to optimize accessibility and travel times in transit networks, namely the Hub Covering Location Problem (HCLP). We developed two mathematical formulations, $F_1$ and $F_2$, and a matheuristic procedure. 
The HCLP is closely related to covering (hub) location and connected facility location problems, and simultaneously considers two types of coverage: network access distance and maximum network transit time. Formulation $F_2$ outperforms $F_1$, yielding optimal solutions for instances up to 50 nodes in less than two hours of CPU time after preprocessing and the addition of valid inequalities. 
Our matheuristic achieves remarkable computational performance, obtaining optimal or near-optimal solutions in less than a minute. Sensitivity analysis shows, on the one hand, that allowing longer travel times in the transit network results in more expansive hub networks with lower density. On the other hand, increasing the coverage radius has the opposite effect.

Future research directions include extensions to multiple lines or transportation modes (multimodal transportation) and the modeling of transfer times between lines.

\section*{Acknowledgments}
The research of C.-A. Domínguez-Bravo and E. Fernández has been partially funded by the Spanish Ministry of Science and Innovation through project PID2023-146643NB-I00 of the State Research Agency (AEI). A. Lüer-Villagra also gratefully acknowledges the support of grants ANID FONDECYT 1261719 and ANID CATLEC CIN250061.

\bibliographystyle{elsarticle-harv}
\bibliography{References20260503}

\appendix

\section{Detail of best-known solutions attained with $F_2$ warm start variant}

\begin{figure}[htbp]
\begin{subfigure}[b]{0.35\textwidth}
\frame{\includegraphics[width=\textwidth]{10CAB_R340_T230_network.pdf}}
\caption{$|V|=10$}
\end{subfigure}
\begin{subfigure}[b]{0.35\textwidth}
\frame{\includegraphics[width=\textwidth]{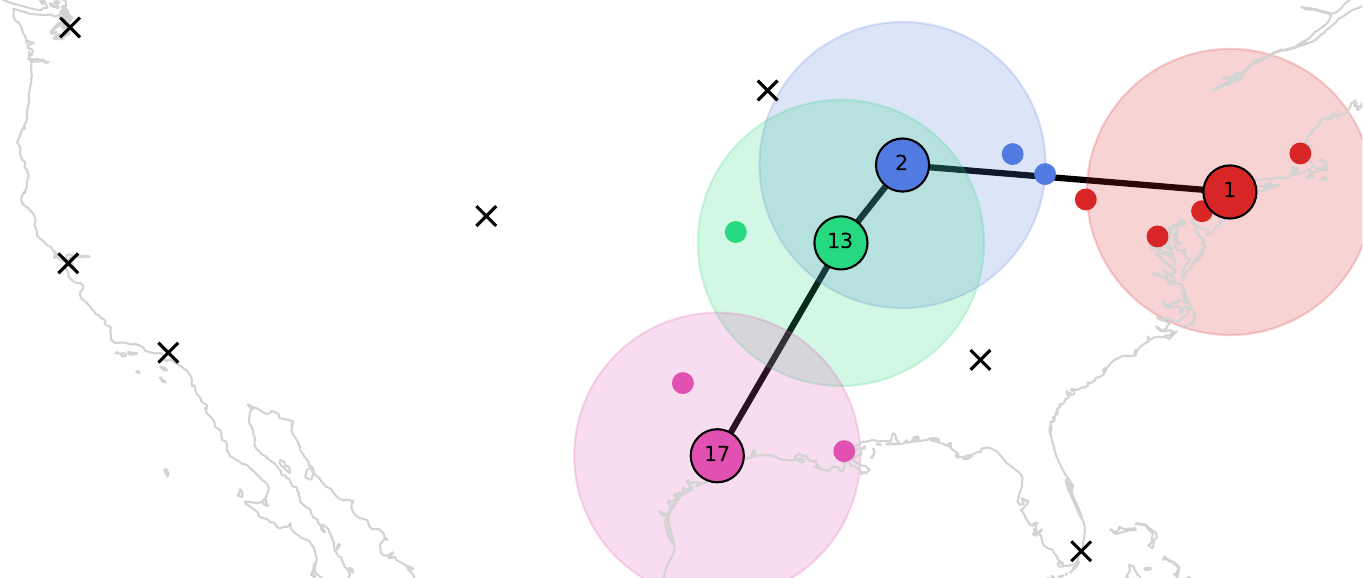}}
\caption{$|V|=20$}
\end{subfigure}
\begin{subfigure}[b]{0.35\textwidth}
\frame{\includegraphics[width=\textwidth]{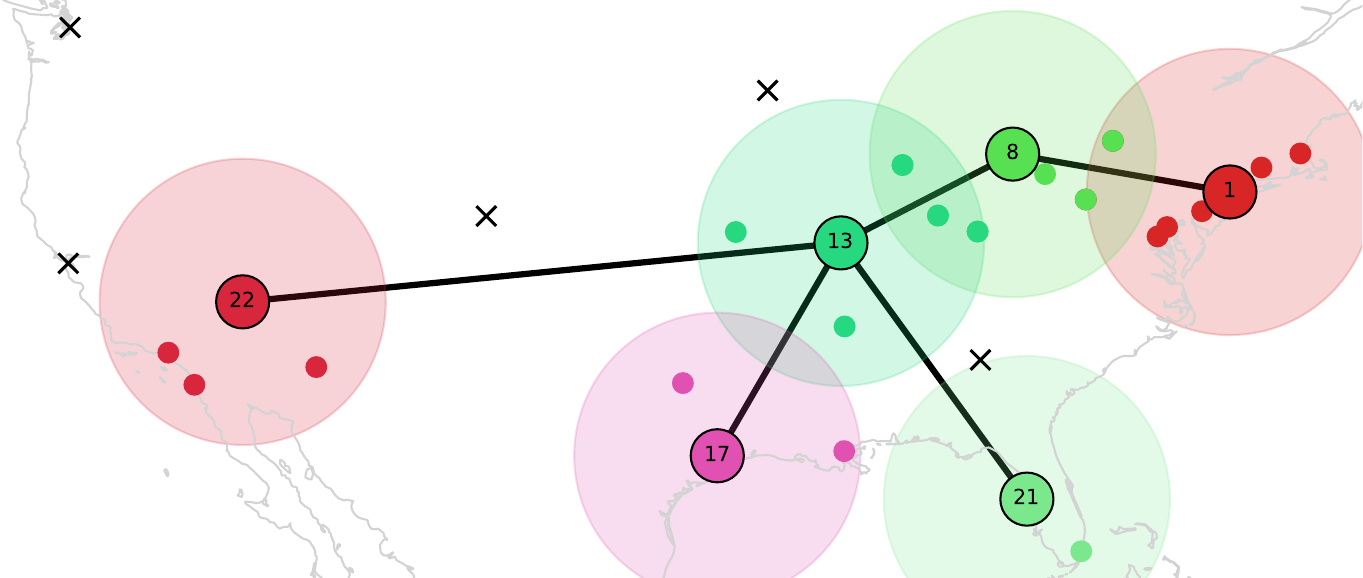}}
\caption{$|V|=30$}
\end{subfigure}
\\
\begin{subfigure}[b]{0.35\textwidth}
\frame{\includegraphics[width=\textwidth]{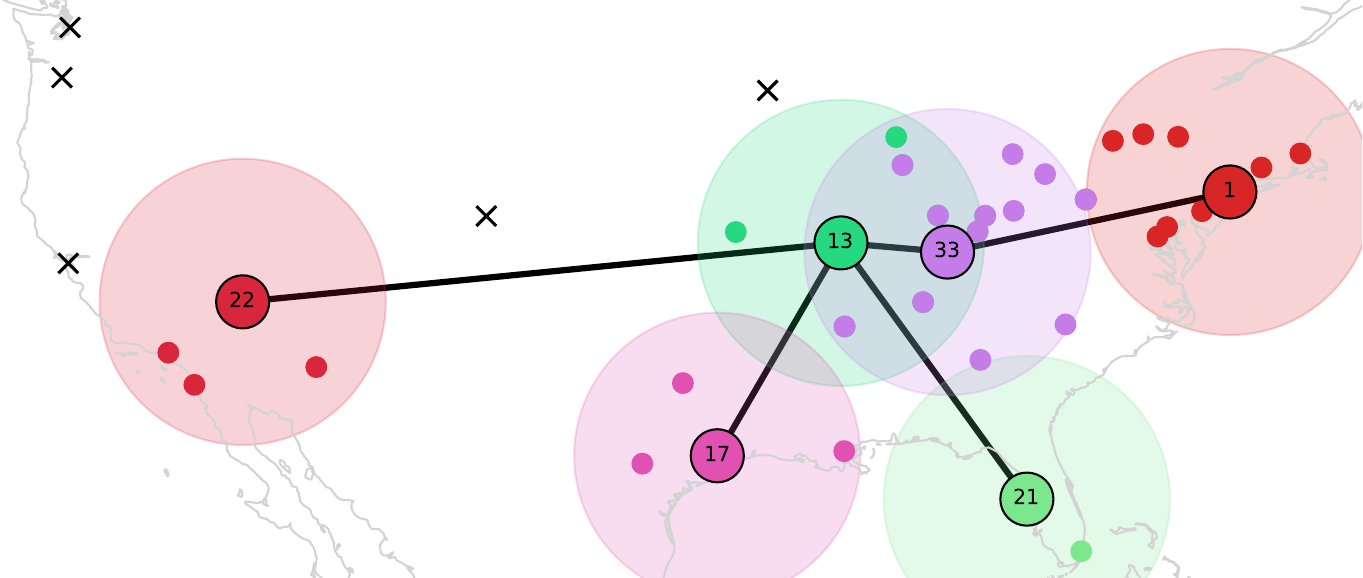}}
\caption{$|V|=40$}
\end{subfigure}
\begin{subfigure}[b]{0.35\textwidth}
\frame{\includegraphics[width=\textwidth]{50CAB_R340_T230_network.pdf}}
\caption{$|V|=50$}
\end{subfigure}
\begin{subfigure}[b]{0.35\textwidth}
\frame{\includegraphics[width=\textwidth]{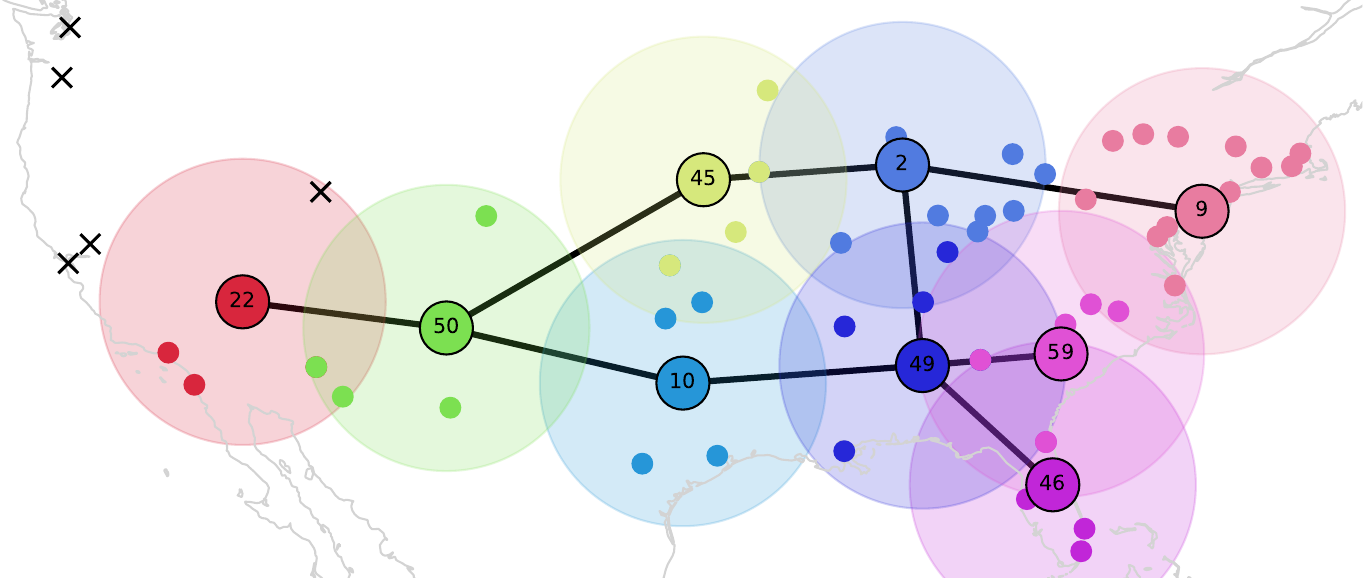}}
\caption{$|V|=60$}
\end{subfigure}
\\
\begin{subfigure}[b]{0.35\textwidth}
\frame{\includegraphics[width=\textwidth]{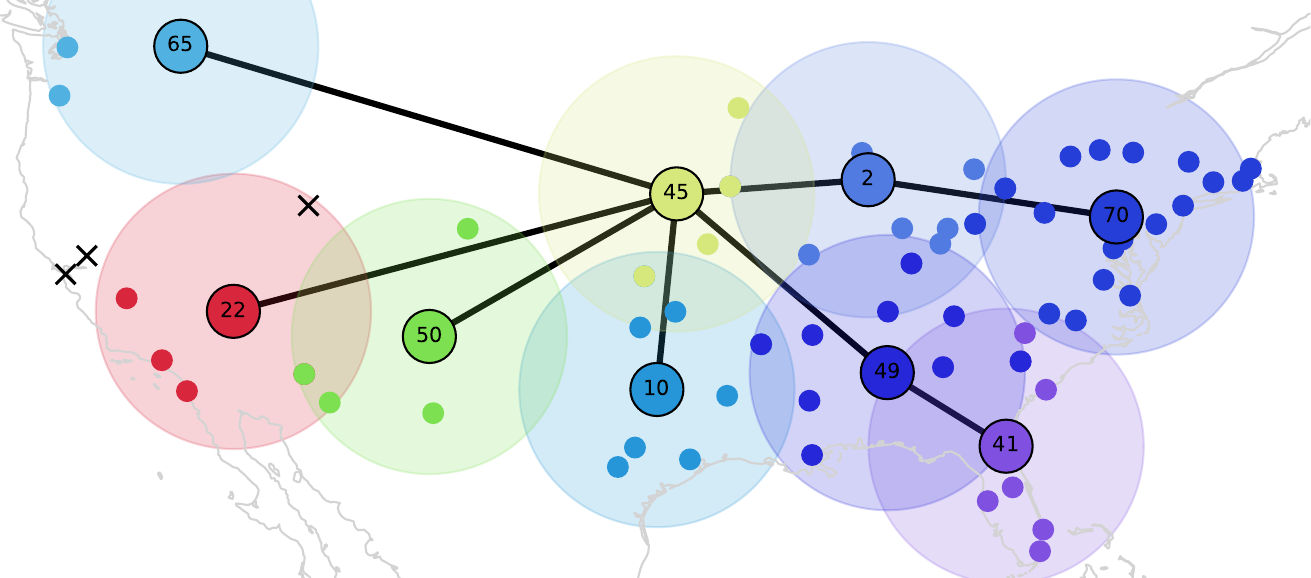}}
\caption{$|V|=70$}
\end{subfigure}
\begin{subfigure}[b]{0.35\textwidth}
\frame{\includegraphics[width=\textwidth]{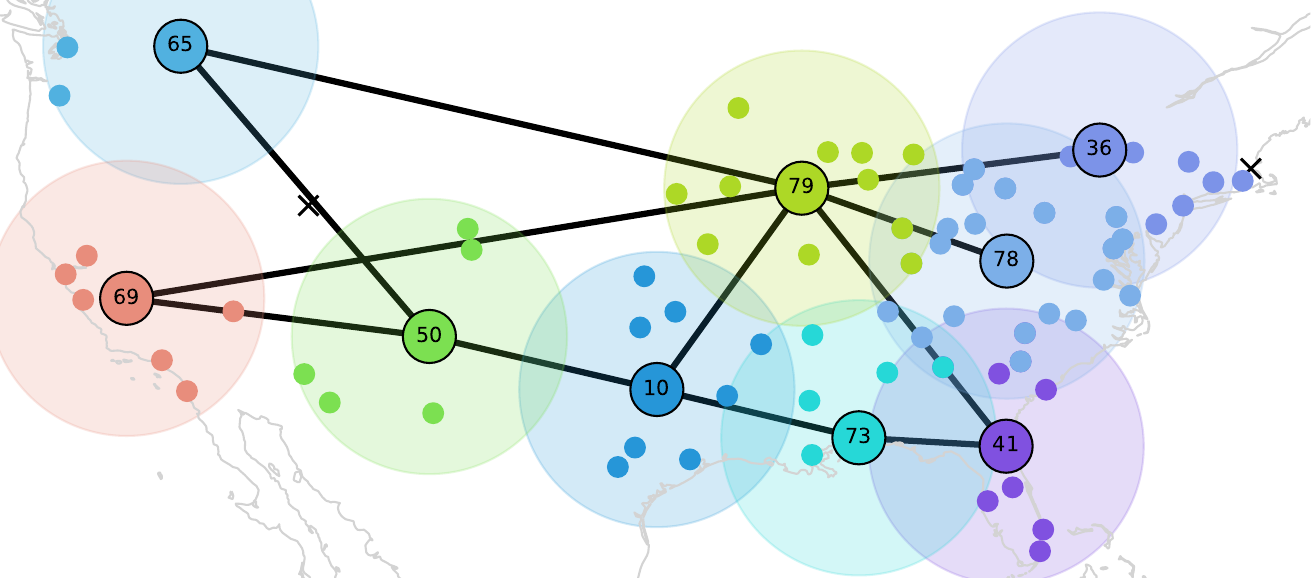}}
\caption{$|V|=80$}
\end{subfigure}
\begin{subfigure}[b]{0.35\textwidth}
\frame{\includegraphics[width=\textwidth]{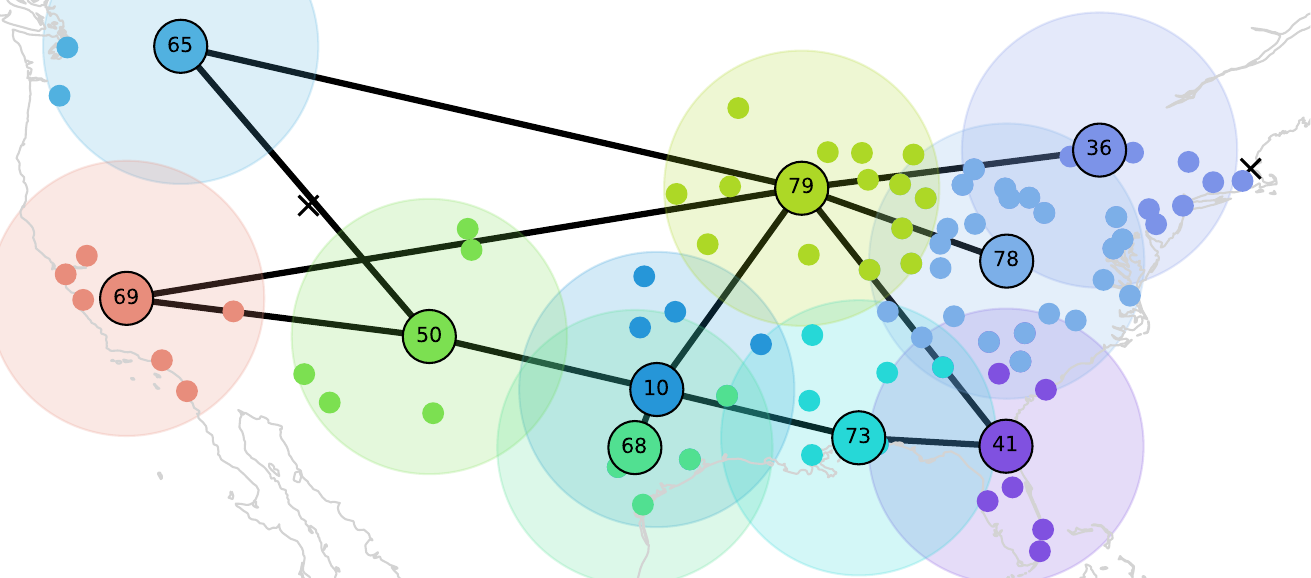}}
\caption{$|V|=90$}
\end{subfigure}
\\
\begin{subfigure}[b]{0.35\textwidth}
\frame{\includegraphics[width=\textwidth]{100CAB_R340_T230_network.pdf}}
\caption{$|V|=100$}
\end{subfigure}
\caption{Best known solutions for CAB instances ($T_{max}=230$ and $\delta=340$).}
\label{fig:best-solutions}
\end{figure}

\end{document}